\documentclass[graybox]{svmult}

\usepackage{mathptmx}
\usepackage{helvet}
\usepackage{courier}

\usepackage{makeidx}
\usepackage{graphicx}
\usepackage{multicol}
\usepackage[bottom]{footmisc}
\usepackage{amsmath,amssymb,url}

\newcommand{\comment}[1]{\relax}

\newcommand{\Hz}{\operatorname{Hz}}
\newcommand{\divergence}{\operatorname{div}}
\newcommand{\Ca}{\operatorname{Ca}}
\newcommand{\Fe}{\operatorname{Fe}}
\newcommand{\nm}{\operatorname{nm}}

\newcommand{\elas}{\operatorname{elas}}
\newcommand{\vis}{\operatorname{vis}}
\newcommand{\ext}{\operatorname{ext}}

\begin{document}

\title*{Chapter 30 -
Geometric and Electromagnetic Aspects of Fusion Pore Making}

\author{D. Apushkinskaya, E. Apushkinsky, B. Boo{\ss}-Bavnbek,
and M. Koch}

\authorrunning{D. Apushkinskaya et al.}

\institute{Darya Apushkinskaya \at Fachbereich Mathematik,
Universit{\"a}t Saarbr{\"u}cken, Postfach 151150, D-66041
Saarbr{\"u}cken, Germany, \email{darya@math.uni-sb.de} \and Evgeny
Apushkinsky \at Experimental Physics Department, St. Petersburg
State Polytechnical University, St. Petersburg, Russia,
\email{apushkinsky@hotmail.com} \and Bernhelm Boo{\ss}-Bavnbek\at
Roskilde University, IMFUFA, Dept. of Science, Systems and
Modelling, DK-4000 Roskilde, Denmark, \email{booss@ruc.dk} \and
Martin Koch\at Feldkraft Ltd., H{\o}jbovej 11, DK-2500 Copenhagen,
Denmark, \email{mko@mail.tele.dk} \at \at File name: 30\_ApApBoKo\_rev,
revised: \today}

\maketitle

\setcounter{minitocdepth}{2}
\dominitoc



\abstract{For regulated exocytosis, we model the morphology and
dynamics of the making of the fusion pore or porosome as a
cup-shaped lipoprotein structure (a {\em dimple} or {\em pit}) on
the cytosol side of the plasma membrane. We describe the forming of
the dimple by a free boundary problem. We discuss the various forces
acting and analyze the magnetic character of the wandering
electromagnetic field wave produced by intracellular spatially
distributed pulsating (and well observed) release and binding of
$\Ca^{2+}$ ions anteceding the bilayer membrane vesicle fusion of
exocytosis. \newline
Our approach explains the energy efficiency of the
dimple forming prior to hemifusion and fusion pore and the
observed flickering in secretion. It provides a frame to relate
characteristic time length of exocytosis to the frequency, amplitude
and direction of propagation of the underlying electromagnetic field
wave.\newline
We sketch a comprehensive experimental program to verify - or falsify -
our mathematical and physical assumptions and conclusions where conclusive evidence
still is missing for pancreatic $\beta$-cells.
 }

\section{Introduction}\label{s:intro}

\subsection{On Our Heuristic (Suggestive) Use of Mathematical Modeling}

This chapter adds a few electromagnetic facts and mathematical
theorems to the toolbox approaching the process of bilayer membrane
vesicle fusion. We address the related geometric and dynamic aspects
of the endocytotic-exocytotic cycle which is at the core of various
discharge (e.g., secretion) and ingestion (e.g., drug intake)
processes in animal cells.

We begin with a {\em caveat}. From mathematical physics, quantum chemistry,
and various fields of engineering design we are accustomed to perfect reliability of
theoretical calculations due to full understanding of the governing laws and full
practical control of the calculated processes. We have
learnt, sometimes the hard way, from physics history that, when in doubt, we had better trust the theory and
carefully designed elaborate experiments than first views and ad-hoc explanations.

Clearly, the situation is different in theoretical biology and medicine. There,
it seems to us, the main use of mathematical modeling is either falsification or extrapolation.
By {\em falsification}, we mean the use of simple
arithmetic or other more advanced mathematical means to check and
falsify common belief (like Harvey's mathematical microscope,
see Ottesen \cite{Ot10}, or the harmonic analysis of $\Ca^{2+}$ oscillations in
Fridlyand and Philipson \cite{FriPhi}).
By {\em extrapolation} we mean the packing of established phenomenology
into a precise, intentionally simplified mathematical frame work admitting
series of computer simulations
or analytical estimates to investigate the r{\^ o}le of selected parameters (like the
Silicon Cell, see Westerhoff et al. \cite{We10}, the mesoscopic simulation of
membrane-associated processes, see Shillcock \cite{Shi}), or the compartment models for
different pools of insulin granules in exocytosis preparations, see Toffolo, Pedersen,
and Cobelli \cite{ToPeCo}).

In this chapter, our approach is different. We simply ask: (1) Could it be that a
highly {\em localizable} phenomenon like the lipid bilayer fusion of regulated exocytosis
on a characteristic length scale of tens of nanometers
has essential {\em cell--global} aspects on a characteristic length scale of hundreds and
thousands of nanometers?  (2) Could it be that the observed
changes of the electrostatic plasma membrane potential accompanying regulated exocytosis
and the corresponding $\Ca^{2+}$ oscillations have
an electro-magnetic character which requires a field-theoretic (Maxwell) approach to
the secretion process?

We have good reason for our, at present still speculative but hopefully suggestive
approach, both in the re-interpretation of more or less well-observed phenomena
and in focusing on aspects which seem to us not sufficiently supported by common
explanations. This will be explained below.

We shall emphasize that the correctness of the electro-dynamical and mathematical modeling
parts of the findings of this chapter depend on future experimental testing and biological
validation. At the end of this report, the reader can find a comprehensive list of experiments
that really need to be done to confirm the relevance of our physical equations and mathematical
modeling. Therefore, our
models do not aim for instant clarification but rather set a scene for alternative
considerations and future observations. That's what we understand by
the {\em heuristic} (suggestive) use of mathematics.

\subsection{Electromagnetic Free Boundary Route to
Fusion Pore Making}

For regulated exocytosis, we model the morphology and dynamics of
the making of the fusion pore or porosome as a cup-shaped
lipoprotein structure (a {\em dimple} or {\em pit}) on the cytosol
side of the plasma membrane. One ingredient to our model is a free
boundary problem for the dimple under the action of electromagnetic
forces, \comment{Formulation OK?} in particular the Lorentz force
acting on charged molecules of the cell's plasma
membrane with decreasing capacitive reactances while forming the
dimple. The force comes from a wandering electromagnetic field
produced by intracellular spatially distributed pulsating (and well
observed) release and binding of $\Ca^{2+}$ ions.

Our approach is based on variational principles and emphasizes
regularity and singularity under the deformation process of the
membranes. It explains the energy efficiency of the observed dimple
forming prior to hemifusion and fusion pore and the observed
flickering in secretion. It provides a frame to relate
characteristic time length of exocytosis (ranging between
milliseconds in nerve cells and seconds in $\beta$-cells) to the
frequency, amplitude and direction of propagation
of the underlying electromagnetic field wave.

We shall not address all the {\em machines} (both protein machines
and lipid assemblies) working together in making the structure and
the composition of the fusion pore. Admittedly, conclusive evidence
is still lacking of the critical character of the here described
electromagnetic field wave for the well-functioning of the regulated
exocytosis in healthy cells and the lack of secretion robustness in
stressed cells. In germ, however, the present electromagnetic free
boundary model gives various hints to future calculations,
estimates, and in {\it vivo}, in {\em vitro} and in {\it silico}
(i.e., numerical simulation) experiments.

\subsection{Plan of the Chapter}
In Sect. \ref{s:synopsis} we summarize several mathematical,
electrodynamical and cell physiological facts which seemingly have
been overlooked or discarded in the literature, but may in our
perception add essential ingredients to a comprehensive
understanding of the short very first phase of regulated exocytosis.
In Sect. \ref{s:model} we describe our model and the corresponding
differential equations, force balances and cost functionals.
In Sect. \ref{s:regularity} we discuss regularity and singularity
results. In Sect.~\ref{s:drive} we present our
preliminary conclusion; some hints regarding the question of what
controls the speed of the process; and a review of experimental
tasks and capabilities to test our hypotheses.

\section{Synopsis of Established Facts}\label{s:synopsis}

We describe the challenge of making the fusion pore; distinguish different mathematical modeling approaches;
and elaborate electromagnetic and geometric phenomena of the very first phase of regulated exocytosis,
namely $\Ca^{2+}$ oscillations, the corresponding slow and low frequent electromagnetic field wave, and the
forming of a dimple in the plasma membrane prior to the bilayer membrane vesicle fusion.

\subsection{Membrane Fusion and the Fusion Pore Challenge}\label{ss:membrane_fusion}
In animal cells, membrane fusion between the plasma membrane and transport vesicles is fundamental for the secretion of macromolecules. In contrast, the opposite event, i.e. the forming of vesicles or endosomes from the plasma membrane is necessary for the uptake of macromolecules and nutrients. The latter process is known as {\it endocytosis}, in which,
in the suggestive words of a renown textbook ``localized regions of the plasma membrane
invaginate and pinch off to form endocytotic vesicles" (Alberts et al., \cite{AlbMol}).
The process of discharge of material, after collecting it in transport vesicles,
is called {\it exocytosis} and is our subject.

\begin{figure}[htb]
\includegraphics[scale=.6]{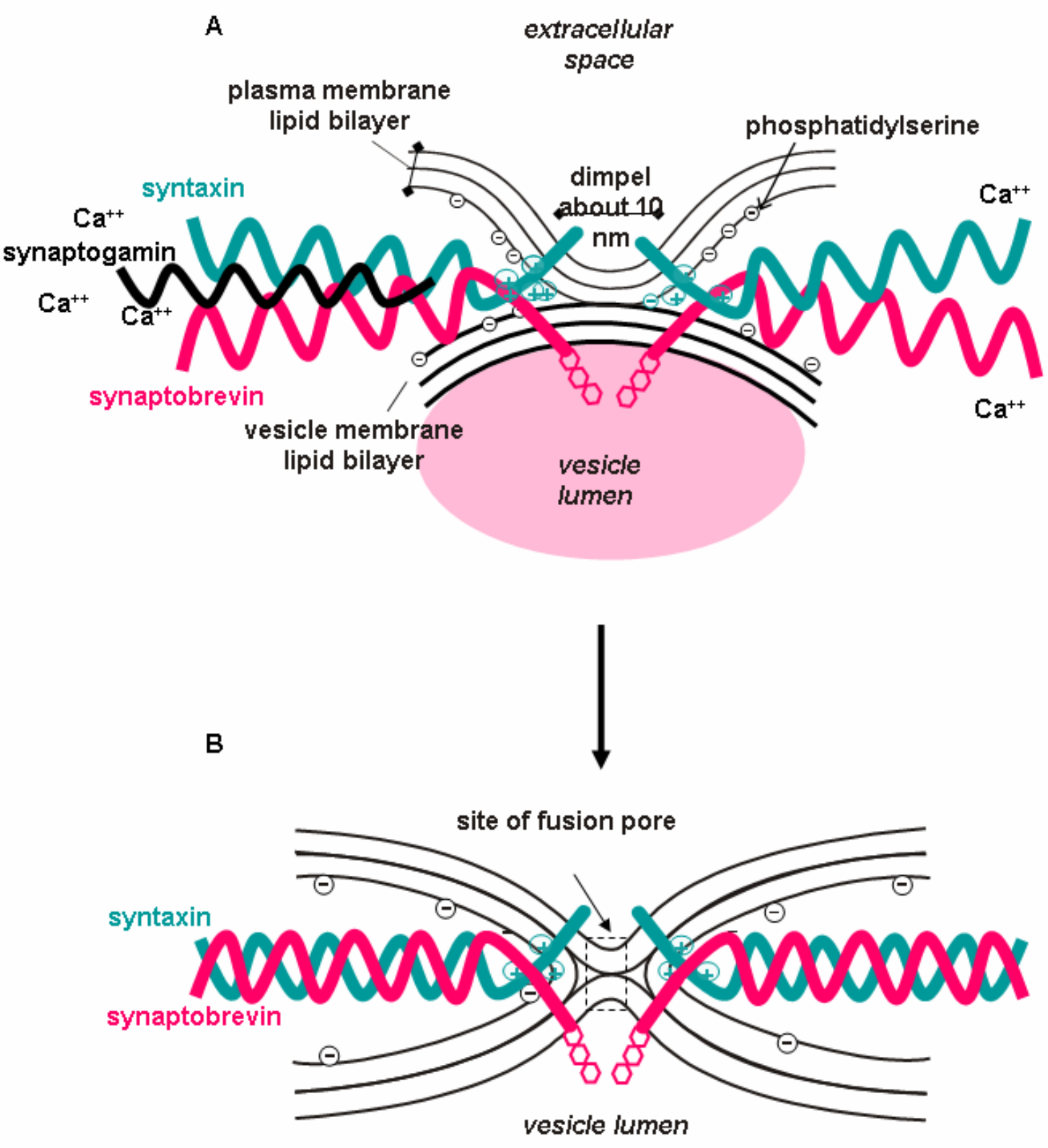}
\caption{Schematic view of the bilayer membrane fusion event, from
Koch et al. \cite{Koc} (after Lentz \cite{Len}), reproduced with
permission of World Scientific Publishers and the authors. Shown are
the basic molecular structures involved in the process of fusion
pore making and membrane fusion. Note the lipid bilayers of the
plasma membrane and the vesicle membrane, with the negative charges
representing the charged heads of phosphatidyl serine. For clarity,
only the transmembrane proteins syntaxin and synaptobrevin with the
three aromatic amino acids on the luminal side of the vesicle
membrane are depicted. Also included is synaptogamin, which binds
$\Ca^{2+}$ and is associated to both syntaxin and synaptobrevin. }
\label{f:exo-scheme}
\end{figure}

It happens by docking of the vesicle to the cell membrane
through activity of several membrane-associated proteins, followed by
vesicle membrane hemi-fusion and the making of a fusion pore in the
membrane through which the material can be released to the exterior,
see Fig. \ref{f:exo-scheme}.
As a rule, the two processes, i.e. endocytosis and exocytosis, seem to be balanced, the one cutting membrane pieces out
of the cell membrane, the other inserting pieces. In both processes, the remarkable is the opening
of the cell {\em without} pinching a hole.

Typical examples of endocytosis are the intake of nutrition, signal molecules, viruses or drugs;
typical examples of exocytosis are secretion of macromolecules such as hormones from endocrine cells, inflammatory mediators from immune cells, neurotransmitter release from nerve cells, or the removal of waste molecules and bi-products. The $\beta$-cells located in the islets of Langerhans in the pancreas are strongly secretory active in producing and releasing insulin.
It seems that a better understanding of these two processes could, e.g., support the early
diagnosis of metabolic diseases, including diabetes mellitus type 2 (exocytic dysfunction, see
Rorsman and Renstr{\"o}m \cite{RorRen:IGD, {RorRen:RIG}}) or a more efficient delivery of insulin analogues
(inducing endocytosis).

Recent advances in observational and manipulative nano techniques and in mesoscopic coarse-grained
computer simulation have provided substantial progress in visualizing,
understanding - and possibly - influencing the bilayer membrane vesicle fusion
(for a recent systematic review we refer to Shillcock and Lipowski~\cite{ShiLipComp}).
With present technology, however, the observational findings are hampered:
geometric shape {\it and} the parameters of change can not be measured to the wanted degree of precision,
{\it simultaneously}. So \cite{ShiLipComp},
p.~S1196 deplores that ``the molecular rearrangements that take place during the final stage of the fusion process,
where the two initially distinct membranes join and produce a fusion pore, cannot yet be resolved by
these experimental techniques" while ``understanding how the stability of lipid membranes is overcome
by the cellular protein machinery when required is a major topic of research" (l.c.).
This is the challenge which we wish to address by summarizing various mathematical, electrodynamical and
empirical facts. Seemingly, these facts have been overlooked or discarded in the literature. We shall show,
however, to what extent these facts yield, in our perception essential, ingredients to a comprehensive understanding
of this short very first phase of regulated exocytosis.

\subsection{Competing Mathematical Approaches to Space-Time Processes}\label{ss:competing}
We shall distinguish three different mathematical-numerical
approaches to modeling spatial-temporal process of regulated
exocytosis: highly aggregated compartment models, spatially
distributed dynamical systems, and space-time integrating partial
differential equations, where our focus will be. Clearly, all three
approaches admit extensions from rigid, stiff and hence fragile
deterministic to more robust stochastic modeling. Here, however, we
shall not discuss such extensions.

A first class of mathematical exocytosis models are {\em
Compartment Models}, first introduced by Grodsky \cite{Gro} in  1972, assuming that there are two
compartments (pools) of insulin granules, docked granules ready for secretion and reserve granules.
By assuming suitable flow rates for outflow from the docked pool and resupply from the reserve pool to the docked
pool, the established biphasic secretion process of healthy $\beta$-cells could be modeled {\em qualitatively} correct.
By extending the number of pools from two to an array of six and properly calibrating all flow rates, Chen, Wang and Sherman in \cite{CWS} obtained
a striking {\em quantitative} coincidence with the observed biphasic process,
see also Toffolo, Pedersen, and Cobelli \cite{ToPeCo}
in this volume. The nice thing by such compartment models is that they invite the experimentalists (both in
imaging and in proteomics) to verify the distinction of all the hypothetical compartments in cell reality and
to assign biophysical values to the until now only tuned flow rates. A self-imposed limitation is the low resolution
of the aggregated compartments which does not allow to investigate the local geometry and the energy balance of the
secretion process.

On the opposite length and time scale, we have a second class of mathematical models, namely
the impressive numerical analysis of the bilayer membrane vesicle fusion
by {\em Molecular Dynamics} (MD), {\em Monte Carlo simulations} (MC) and {\em Dissipative Particle Dynamics} (DPD)
on nanometer distances and fractions of nanoseconds, based on gravitational and electric forces between the particles,
see the afore cited \cite{ShiLipComp} and Shillcock, \cite{Shi} in this volume.
Unfortunately, these computer simulations are also seriously hampered,
namely by limitations of present hardware and software when one is addressing mesoscopic behavior, i.e., changes across many scales of the molecular characteristics - in spite of the impressive results when applying these methods to phenomena on the nano scale, like modeling the island dynamics of film crystallization in epitaxial growth driven by molecular beam epitaxy, see, e.g., Caflisch and Li \cite{CafLiAnal}. These limitations in present computer capacity require MD, MC
and DPD simulations to make {\it a priori} assumptions about the {\it pathway} of the fusion process, e.g., spherical symmetry of the vesicles and planarity and circularity of the fusion pores - besides the often deplored
``enormous gap between the sophistication of the models and the success of the numerical approaches used in practice and, on the other hand, the state of the art of their rigorous understanding" (Le Bris \cite{LebICM} in his 2006 report to the International Congress of Mathematicians). To keep these models transparent, a self-imposed limitation is their
focus on the local neighborhood of the fusion event, neglecting long-distance phenomena like electromagnetic waves across the cell.

Whereas compartment models respectively MD and DPD are built upon small respectively huge
systems of ordinary differential equations with each unknown specifying temporal changes in one given
pool or spatial box, we advocate a third class of mathematical models, namely modeling the dynamics and the geometry
by {\it partial differential equations}. Consequently, we shall try to model
the relevant processes by one or two spatial-temporal equations based on First Principles,
instead of the few aggregated purely temporal pool equations in compartment modeling or the three millions of
purely temporal equations for spatially distributed boxes
in \cite{ShiLipComp}, p. 1197 (which still is a very poor particle number for a 100 nm $\times$ 100 nm $\times$ 42 nm simulation box).

Moreover, the simplicity of our fundamental equations admits a transparent incorporation of long-distance
phenomena. In such a way, our approach takes its point of departure not only in the rather well studied
{\em elastic} and {\em electric} properties and potentials of and across the plasma membrane
and the {\em viscosity} of the cytosol,
but in non-stationary, dynamic {\em electromagnetic} properties.
To us, the basic electromagnetic character of the fusion process becomes evident in
\begin{itemize}
\item the observed electromagnetic (wandering) field waves, see Sect. \ref{ss:field-wave},
\item the closure of the corresponding magnetic wave over the plasma membrane, see Sect. \ref{ss:magnetic-closure},
\item the observed forming of a narrow dimple, solely to be understood like a capacitator, see
Sect. \ref{ss:dimple}, and
\item the observed flickering of the secretion process corresponding
to the natural variability of the (wandering) field wave generation, see Sect. \ref{ss:flickering}.
\end{itemize}
The goal of our approach is
\begin{enumerate}

\item to develop a simple free boundary model for the dimple forming process, see Sect. \ref{s:model},

\item to focus on regularity and possible singularity of the free
boundary, see Sect. \ref{s:regularity},

\item to provide a reliable framework for estimating (and, hopefully, influencing)
the parameters which control the speed of the process, see Sect. \ref{ss:findings},

\item and to formulate a bundle of model-based observation plans to verify or falsify our assumptions,
see Sect. \ref{ss:suggested}.

\end{enumerate}

Our approach is inspired by recent work of A. Friedman and collaborators about tumor growth,
see \cite{BazFriCPDE}--\cite{TaoCheN};
by a theoretical analysis of adhering lipid vesicles with free edges,
see \cite{NiShiYinColl} by Ni, Shi, and Yin; by the electrodynamic challenge to understand
the observed $\Ca^{2+}$ oscillations, rightly perceived as being ``contradictory and often
do not support the existing (electrostatic) models" (Fridlyand and Philipson
\cite{FriPhi}) and
by the mathematical challenge to understand the (in elastic terms counter-intuitive)
dimple formation so well described in the literature,
see, e.g., Monck and Fernandez \cite{MonFer}, Rosenheck \cite{Ros},
Lentz et al. \cite{Len}, Koch et al. \cite{Koc}.

\subsection{Oscillatory Intracellular Release and Binding of $\Ca^{2+}$
Ions}\label{ss:field-wave}
We recall a few basic observations of $\Ca^{2+}$-oscillations and
postulate a simple but powerful procedure of generating an electrodynamic field.

\subsubsection{Basic Observations of $\Ca^{2+}$-oscillations}
By fluorescence microscopy, empirical evidence has been provided
about pulsating $\Ca^{2+}$ activity at extreme low frequency $f\sim
0.1 \Hz \ll 3 \Hz$
(for comparison, the house low frequency grid is
of $50 \Hz$, i.e., spikes in intervals of 20 msec), prior to the
fusion event, see Kraus, Wolf and Wolf, \cite{KrWoWo} and Bernhard
Wolf's homepage \cite{Wo} with informative video animations of
calcium oscillations, and the comprehensive review and analysis by
Salazar, Politi and H{\"o}fer, \cite{SPH}.\footnote{Our model cell
is a pancreatic $\beta$-cell where a single release is slow and may
take seconds. Correspondingly, we expect a low $\Ca^{2+}$
oscillation rate in intervals in the range of seconds yielding
extreme low frequency of the observed $0.1 \Hz$. For nerve cells the
reaction time, and so the release time, is in the range of msec,
possibly less than 100 $\umu$sec, see Jahn, Lang, and S{\"u}dhoff
\cite{JLS}. Correspondingly, we expect a high $\Ca^{2+}$ firing rate
in intervals of, e.g., 10 msec yielding a frequency of 100 Hz with
associated high energy losses. So, our electromagnetic free boundary
route to vesicle fusion can not function in nerve cells unless the
neurotransmitter vesicles are kept waiting very close to the plasma
membrane. Here our finding coincides with the well-known deviating
high energy consumption of nerve cells. We shall elaborate this
aspect in Sect. \ref{s:drive} below.}
Wolf's observations were made with HELA cancer cells. H{\"o}fer et al.
gives a general model, based mainly on observations in muscle cells.

For corresponding observations for $\beta$-cells, we refer to Maechler
\cite{Ma10} and Fridlyand and Philipson \cite{FriPhi} in this volume,
who deal with
various types of low-frequent oscillations.

The following can be seen in many cell types:
when a cell is stimulated (``polluted") by a $\Ca^{2+}$--mobilizing stimulus, the
changes in the cytosolic calcium concentration $[\Ca^{2+}]_c$ occur
as repetitive spikes that increase
their frequency with the strength of the stimulus (see also Berridge, Bootman,
and Roderick \cite{BBR}
and Gaspers and Thomas \cite{GasTho}). It is well known that an increase in
$[\Ca^{2+}]_c$\,, ultimately, regulates a plethora of cellular processes mediated by
$\Ca^{2+}$-dependent enzymes that, in turn, modify downstream
targets commonly by phosphorylation.

Investigating $\Ca^{2+}$ decoding in an analytically tractable generally applicable
model, H{\"o}fer and collaborators address the question ``Under
which conditions are $\Ca^{2+}$  oscillations more potent than a
constant signal in activating a target protein?" in \cite[p.
1204]{SPH} in complex biochemical terms of binding and release rates
of $\Ca^{2+}$  ions.

\subsubsection{Postulated Electrodynamic Field Character}

We wish to supplement these investigations by recalling
pieces of circumstantial evidence which may support
the hypothetical electromagnetic character of the $\Ca^{2+}$
oscillations. As said before, $\Ca^{2+}$ handling is extremely complex.

\begin{enumerate}
\item There are various $\Ca^{2+}$-storage organelles in the $\beta$-cells, first of all
the spatially distributed and clearly separated mitochondria and the spatially rather extended
smooth endoplasmic reticulum - SER, see also \cite[Section
5]{FriPhi}.

item There is some agreement in the literature that the  handling
by the mitochondria, and not the SER, is decisive for, at least, the slow
$\Ca^{2+}$ oscillations: ``Metabolic profiling of $\beta$-cell function identified
mitochondria as sensors and generators of metabolic signals controlling insulin secretion",
according to Maechler \cite{Ma10}; and Fridlyand and Philipson \cite{FriPhi}
refer to ``data demonstrating that slow oscillations can persist in the presence of
thapsigargin, the agent that blæocks SERCA and empties the ER stores...".

\item To us, these observations fit nicely with our postulate of the electrodynamic
character of the slow $\Ca^{2+}$ bursts: spatially distributed and temporarily
coordinated (not necessarily simultaneous) release and uptake of $\Ca^{2+}$ ions can
generate an electrodynamic field with magnetic character.

\end{enumerate}

Of course, there are many mechanisms and systems that contribute to the release and
uptake of $\Ca^{2+}$ from intracellular organelles, some of these even have bell-shaped
effects on $\Ca^{2+}$ release, and phenomena such as $\Ca^{2+}$-induced $\Ca^{2+}$
release also take place (these facts were communicated to the authors by Pociot and
St{\o}rling \cite{PoSt}). In the following, we anticipate
the existence of a (not yet fully confirmed) system of aggregated build-up of an
electrodynamic field with magnetic character in $\beta$-cells upon glucose stimulus: In short,
we have movements of $\Ca^{2+}$ ions in and out of the mitochondria. Movements of ions
are currents. In-and-out movements are AC currents. AC currents produce AC fields,
or rather, AC field oscillations. Superposed oscillations may produce moving fields.

In this way, we suppose, the organelles build an
alternating electric current density (also ``displacement vector
field") $\mathbf D$ of low
frequency by superposition of spatially distributed, temporally
coordinated and directed $\Ca^{2+}$ activity. As usual in
electrodynamics, we shall speak of two different electric fields,
$\mathbf D$ and $\mathbf E$. This second electric field $\mathbf E$
is given by the relation $\mathbf D = \varepsilon \mathbf E$, where
$\varepsilon = \varepsilon_0 \varepsilon_r$ denotes the dielectric
constant.  Note that our writing of all electromagnetic units
and equations follows Jackson \cite{Jac} in the units V, A, sec, and
m with, e.g., kg=VAsec$^3$/m$^2$ of the System International - SI,
which is predominant in engineering literature.

A physical model of generating electromagnetic field waves by
spatially and temporally distributed excitation was built by Koch
and Stetter, see \url{http://www.feldkraft.de/}. It is called {\em
Dynamical Marker} and consists of a couple of coil arrays,
electronically regulated for direction-, frequency- and
amplitude-controlled generation of a field wave. The instrument has been
applied in various cell physiological experiments for slow and
efficient transport of beads into cells across the plasma membrane.

There are a couple of related questions which will require separate
investigation (see also below Sect. \ref{ss:suggested}):

\begin{enumerate}

\item How strong is the evidence that the $\Ca^{2+}$  {\em oscillations}, contrary to the
membrane process of $\Ca^{2+}$ bursts,  originate from an
array of $\Ca^{2+}$  depots (SER and mitochondria), organized in directed lines?

\item How does the cell select the $\Ca^{2+}$  storages to participate in the generation of the
alternating current?

\item How is the sequential release and binding of the $\Ca^{2+}$  ions of the different
storages controlled, i.e., how is the correct spatial and temporal coordination of
release and binding obtained?

\item What role plays the observed branching of mitochondria in
active $\beta$-cells,
contrary to the dipole shape of mitochondria in tired $\beta$-cells?

\begin{figure}[htb]
\includegraphics[height=18cm,width=14cm]{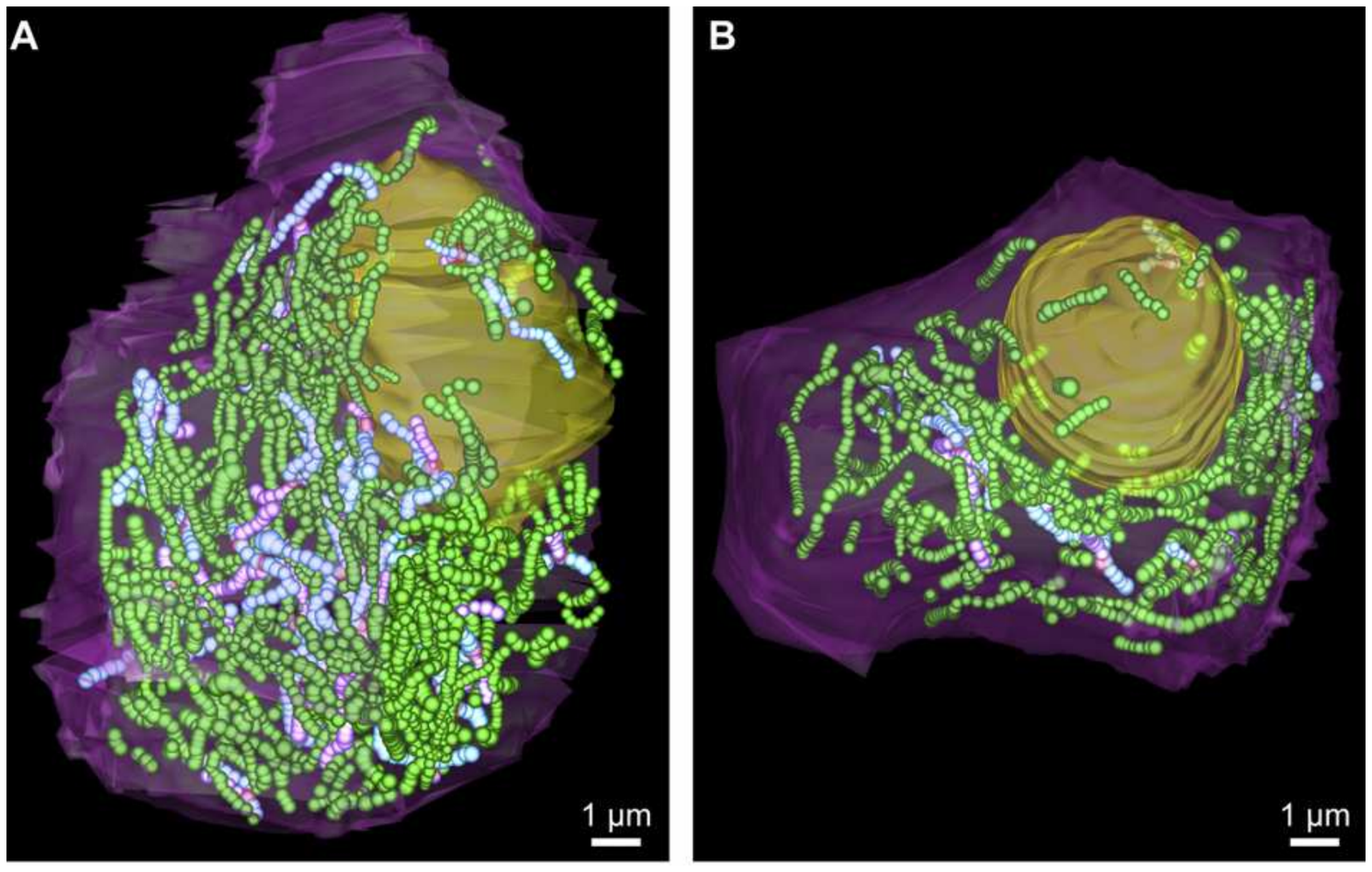}
\vspace{-11.7cm} \caption{Functional heterogeneity: left, a
multitude of {\em branched} mitochondria in a vigorously responding
$\beta$-cell; right, relatively more non-branched {\em dipole}
mitochondria in a less active cell at a comparable locus; branch
points are highlighted with red spheres. The panels are from Noske
et al., \cite{NCMM} and reproduced with permission of Elsevier. High
resolution originals were supplied by courtesy of B. Marsh,
University of Queensland, Brisbane, Australia.}
\label{f:mitochondria}
\end{figure}

\item Can the magnetic character of the field wave, which is
produced by low frequent $\Ca^{2+}$  oscillations, be influenced via
an external field with similar character?

\end{enumerate}

To Questions 1 and 2, it may be mandatory to distinguish between two
different types of $\Ca^{2+}$ burst: the $\Ca^{2+}$  {\em
oscillations} addressed here are prior to the secretion and are,
presumably, generated by arrays of activated calcium depots
distributed through the full length of the cell. On the contrary,
the $\Ca^{2+}$  {\em influx} through the ion channels of the plasma
membrane during regulated exocytosis is mainly effective {\em close} to the
plasma membrane where it increases the concentration $[\Ca^{2+}]_c$
and changes the electric potential across the plasma membrane. It is
worth mentioning that the array character is evident from the
observed oscillations in the form of {\em directed, oriented} waves.

To Question 3, we imagine that the coordination of the activity of the
participating $\Ca^{2+}$  storages is not controlled externally,
e.g., by the nucleus, but happens spontaneously by selforganization: we notice that the
storages of molecular calcium sense and respond to stimuli by periodic release and binding of $\Ca^{2+}$  ions and
suppose that the sequential coordination between spatially distributed loci of release and binding
minimizes the energy consumption for maintaining the activity and
for establishing a suitable average concentration $[\Ca^{2+}]_c$ \,.

To Question 4, we recall from Noske et al \cite{NCMM} a remarkable discovery.
The authors imaged and reconstructed two $\beta$-cells from the same glucose-stimulated mouse islet by single axis,
serial section electron microscope tomography (ET) at magnifications of $4700\times$ and $3900\times$, respectively, that
resulted in whole cell tomograms with a final resolution of 15-20 nm. In addition, they developed several new methods
for the abbreviated segmentation of both cells' full complement of mitochondria (i.e., the most prominent
$\Ca$ storages)
and insulin secretory granules for comparative analysis. 3D reconstruction by ET of each of the two $\beta$-cells
(designated 'ribbon01' and 'ribbon02') indicated that ribbon01 responded more vigorously to glucose-stimulation
than ribbon02 and contained about twice as many branched mitochondria (26 out of a total number of 249 mitochondria)
as ribbon02 (with 10 branched mitochondria out of a total of 168 mitochondria). See also the recent
Marsh and Noske \cite{MarNos} in this volume.

Now, from our electromagnetic point of view, the advantage of branched
mitochondria is clear for the generating of a field wave (i.e., the
pulsating $\Ca^{2+}$ oscillations through the whole length of the
cell): a single (non rotating) dipole can not generate or initiate a
(directed, oriented) field wave. That {\em requires} a branched
structure with spatially and temporally shifted serial activity, as demonstrated
also in the design of the mentioned simple field generator by Koch and Stetter.
Of course, nature's regulation of the ion firing may be much more sophisticated
than the crude engineering design of the ``Dynamic Marker". After all, the eukaryotic cells
had many more years to test and optimize different designs in evolution.

To Question 5, we refer to an experimental setting described below in Sect. \ref{ss:suggested}.

\subsection{The Magnetic Character of the Induced Field Wave}\label{ss:magnetic-closure}

Perhaps one of the most delicate of Maxwell's equations (cf. Box) is
his modification of Amp{\`e}re's Law by adding the displacement
current density $\frac{\partial\mathbf D}{\partial t}$ on the right,
i.e., the electric side:
\begin{equation}\label{e:max}
\operatorname{curl}\, {\mathbf H} = \mathbf{J} + \frac{\mathbf
D}{\partial t}, \text{ respectively } \int_C \mathbf H\, d\mathbf{s}
= \int\int_A (\mathbf J + \frac{\partial \mathbf D} {\partial t})
d\mathbf{A}.
\end{equation}
Here $\mathbf H$ denotes the oriented, electrically generated
magnetic field, $\mathbf B=\mu_0\mu_r \mathbf H$ the
corresponding magnetic flux density (also ``magnetic induction") of
frequency $f$ and amplitude $\hat {\mathbf B}$, and $\mathbf{J}$ the
current density vector for a conductor crossing an area $A$ which is
bounded by a contour $C$.
Moreover $\mu_0$ and $\mu_r$ denote the
absolute and relative magnetic permeability.

 \begin{figure}[htb]
\includegraphics[scale=.68]{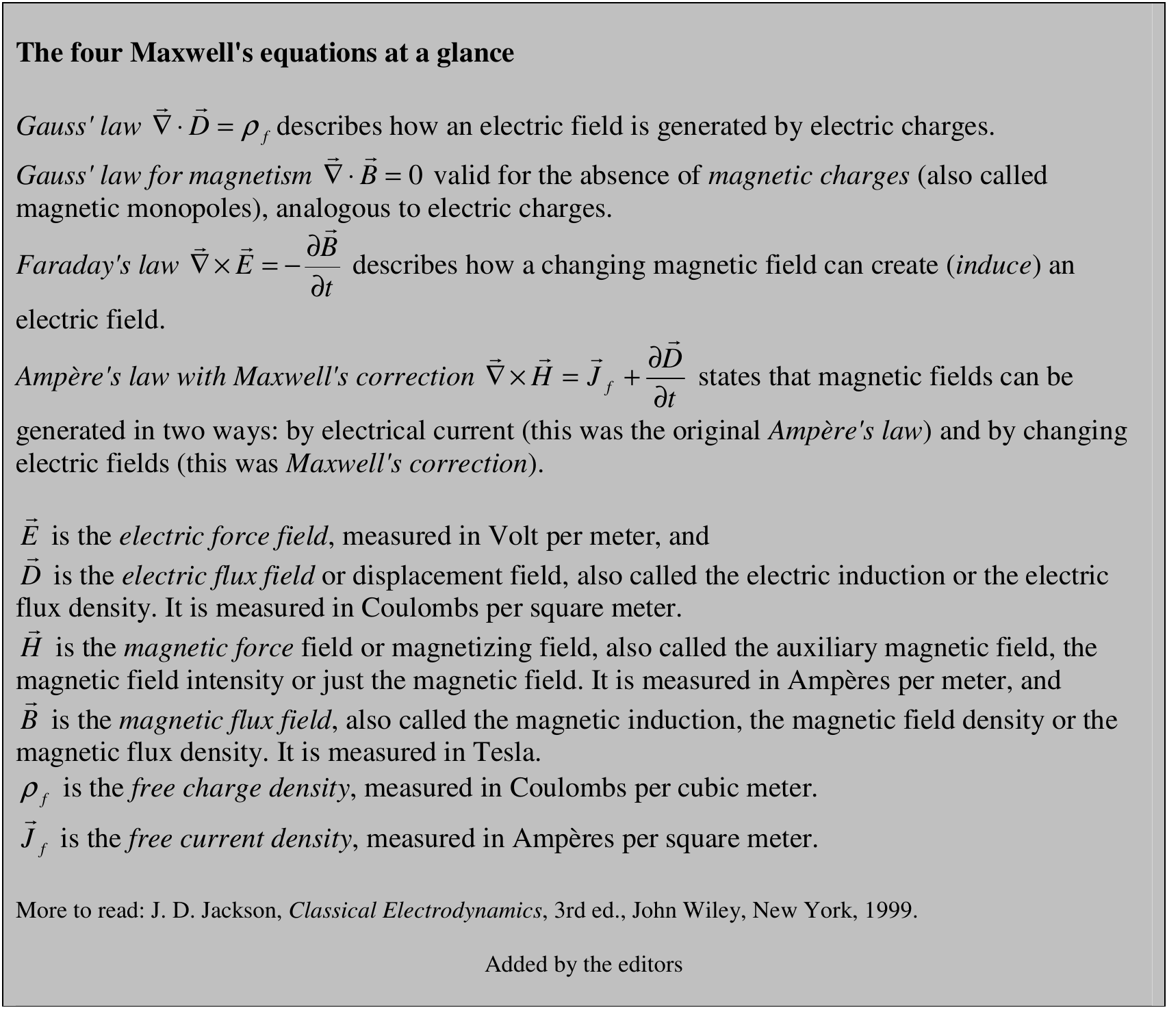}
\end{figure}

Equation \eqref{e:max} is our basic equation regarding the magnetic
character of the observed $\Ca^{2+}$ oscillations. Note that our
field wave can not be compared with an electromagnetic high
frequency wave in radio transmission. Its propagation velocity is
comparable to sound waves in water and far below light velocity.
Moreover, because of the low frequency of our oscillations, the
displacement current density $\frac{\partial\mathbf D}{\partial t}$
is relatively small. That indicates that the magnetic character
dominates the electric character of the field wave.

While a wave with dominant electric character  has large losses in
cytosol (which is comparable to salty water in its electric
conductivity), the magnetic character makes the propagation almost
free of losses. To sum up, by the described intracellular $\Ca^{2+}$
oscillations a field is generated with only marginal losses because
the transmission is almost independent of the material constant
$\varepsilon_r$\,. As a result, the moving $\Ca^{2+}$ ions from the
intracellular distributed storages are in fact an AC current,
generating a nearly loss-free moving (wandering) magnetic {\em field
wave} which transfers energy to a selected transmembrane subregion
of the cytosol between the plasma membrane and a single vesicle. As we
shall explain, this energy and the corresponding forces act on the
free ions distributed in one or other way among the phospholipids of the plasma membrane
and pull the plasma membrane
towards the vesicle. See also \cite{Ko} for a video animation of a
macroscopic field wave.

Our biophysical approach is classical, hence we assume that there
are no sources for magnetism, no magnetic monopoles, at least not
present in our $\beta$-cell. Consequently, we obtain from Maxwell's
equations $\divergence\mathbf B=0$, i.e., the magnetic field wave
induced by the alternating current is {\em closed} in the sense of
vector analysis and, consequently, the path of the field wave is
closed.

Of course, it must be investigated in detail {\em how} the magnetic wave is closed. Fluorescence microscopy
gives the impression that the
observed $\Ca^{2+}$  oscillations are collective phenomena of cell ensembles: the oscillations propagate
through the ensemble like chained waves. We, however, assume that all magnetic waves are separated from
each other and are closed over the single plasma membranes. One reason is that the plasma membrane is perforated by a multitude of ion channels created and maintained by the presence of enzymes like
various kinases and phospholipases. Most encymes contain $\Fe$ atoms, see Jensen \cite[pp. 134f]{Jen}
and the review on iron biominerals \cite{FraBla}.
Consequently, the magnetic field wave will search for a circuit through the plasma membrane.

While we have put the determination of the $\Fe$ content of the
plasma membrane on our experimental
agenda below in Sect. \ref{ss:suggested}, we should mention that the magnetic
permeability $\mu_r$ can not be
measured directly. However, there are methods to determine the magnetic
permeability $\mu_r \sim 1.0000007$ of hemoglobin of deoxidized venous blood non invasively
in a precise way. Similar methods will be applicable for investigating
the plasma membrane of life cells. Roughly speaking, the indirect methods work by comparative measurements after
inflicting a magnetic pollution on a harmonic oscillator.

Now it is not difficult to understand the making of the fusion pore and the dimple formation (see
the following Sect. \ref{ss:dimple} for the empirical evidence) in qualitative terms. Let us fix the
notation. By electron microscopy, we can distinguish the following clearly separated regions:
  \begin{description}
  \item[$D_0$] Amorphous outside cell neighborhood
  \item[$D_1$] Plasma membrane with boundary $\partial D_1=\Gamma_{1-0}\cup\Gamma_{1-2}$
  \item[$D_2$] Cytosol
  \item[$D_3$] Vesicle membrane with boundary $\partial D_3=\Gamma_{3-2}\cup\Gamma_{3-4}$
  \item[$D_4$] Vesicle lumen
  \item[$M_j$] Activated molecular $\Ca$ storage organelles, $j=1,\dots n$ (not depicted)
  \item [$N$] Cell nucleus (not depicted).
  \end{description}

\begin{figure}[htb]
\sidecaption[t]
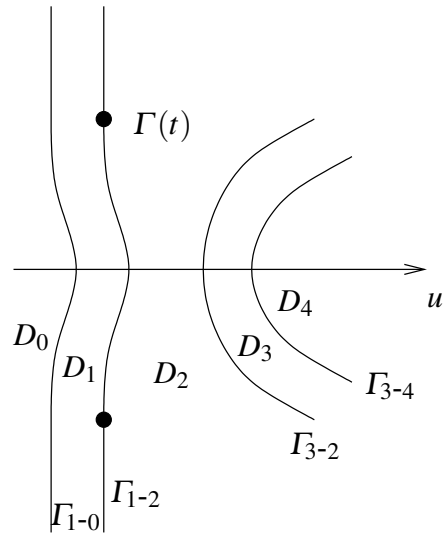
\caption{Plasma membrane section with cone like dimple (left) under
transversal displacement $u$ and spherical vesicle (right) before
docking. For a regular cone shaped dimple, the two bold dots mark
the circular base line $\Gamma(t)=\partial\{(x,y,z)\mid
u(x,y,z,t)\ge 0\}$ of the dimple at time $t$. In space-time, the
union $\{\Gamma(t)\times \{t\}\}_{t\ge 0}$ of the $t$-components
$\Gamma(t)$ form the free boundary $\Gamma$, see below Sect.
\ref{ss:free-boundary} before Eq. \eqref{sign-of-ut} and Sect.
\ref{ss:class_gamma}} \label{f:vesicle}
\end{figure}

\comment{relabel axis with $z$ instead of $x_1$} Note that the
region $D_1$ (the plasma membrane) consists of a phospholipid
bilayer. It has a surface $\Gamma_{1-0}$ as its outside boundary
(towards the amorphous outside cell neighborhood $D_0$) and a
surface $\Gamma_{1-2}$ as its inside boundary (towards region $D_2$
of the watery cytosol within the cell close to the plasma membrane).
Moreover, we have the region $D_3$ consisting of the vesicle
membrane (same material like $D_1$) and the region $D_4$, the
interior of the vesicle, containing the material to be released
through the plasma membrane. Moreover, there is a multitude of
activatable $\Ca^{2+}$ storage organelles $\{M_j\}$ spread through
the interior of the cell. Finally, we have the cell nucleus $N$, see
the abstraction of Fig. \ref{f:exo-scheme} in Fig. \ref{f:vesicle}.
Typical approximate diameters are 10 $\mu m$ for most animal cells;
5 $\mu m$ for the nucleus; 100 $\nm$ for the vesicles; and 8 $\nm$ for
plasma and vesicle membrane.

In the preparation of the fusion event and the making of the fusion
pore, there are apparently only two active regions, in addition to
the cell nucleus and the mitochondria and other $\Ca$ storages,
namely the plasma membrane $D_1$ and the cytosol $D_2$. The plasma
membrane $D_1$ forms a conical inside oriented {\em dimple} (pit)
towards the vesicle of around 10 $\nm$ base diameter and 10-20 $\nm$
height. For the true lipid bilayer membrane-vesicle fusion event,
following the making of the dimple, so-called transmembrane proteins
become active in the cytosol region $D_2$ between plasma membrane
dimple and the vesicle and pull and dock the vesicle membrane $D_3$
to the plasma membrane $D_1$ over a distance of up to 100 $\nm$ (see
also Fig. \ref{f:exo-scheme}).

When the deformation of the plasma membrane is sufficiently sharp it
ends in a branch point, i.e., a singularity of $\Gamma_{1-2}$,
called hemifusion. Then it comes to a break-through (called fusion
pore), and the content of the vesicle begins to diffuse from the
vesicle compartment $D_4$ into the outside region
$D_0$. It appears that this process sometimes is interrupted (so
called {\em flickering}, see below), i.e., the fusion pore is hardly
maintained by elastic forces alone but needs probably the presence
of an electromagnetic field and is interrupted when this magnetic
field is interrupted.

What is controlling the well-functioning of the fusion event?

{\bf Working hypothesis 1}: The regions
\begin{equation}\label{e:reg}
D_0 \text{ ($X_C=0$),\quad
$D_1$ (variable $X_C(x,t)$) \quad and\quad $D_2$ ($X_C=0$)}
\end{equation}
are distinguished by their {\em capacitive
 reactance} $X_C := 1/(\omega C)$ where $\omega = 2\pi f$ with constant $f$ and $C$
 denotes the capacitance. Note that forming the dimple produces an increase of the dielectricum
 (between the ``plates") and so implies increasing $C$
and decreasing $X_C$ until $X_C$ vanishes in the fusion pore.

{\bf Working hypothesis 2}:
The vesicle is densely packed with material and so not subjected to deformations easily.

{\bf  Working hypothesis 3}: We envisage the following feedback mechanism for forming the dimple and preparing the fusion event.
\begin{enumerate}
\item   The $\Ca^{2+}$ ions from locally distributed intracellular $\Ca$-molecule storages start low frequent oscillations,
as described in \cite{KrWoWo}, which are superposed in a controlled way, and
a dynamic field wave is produced pointing to a specific region
$D_{2,\operatorname{crit}}$ selected for most suitable membrane-vesicle fusion.
In the beginning, the magnetic  flux
density vector $\mathbf B$ (``magnetic induction") is low because
the magnetic wave does not easily enter the plasma membrane $D_1$ to
close itself in a circuit because of the originally high $X_C$ in
$D_1$.

\item The form change decreases $X_C$ close to the emerging dimple. That permits the magnetic wandering wave to enter $D_1$ more easily and so increases its current density (the sharpness of its pointing) and its amplitude $\hat B$. And so on.

\end{enumerate}

\begin{figure}[htb]
\sidecaption[t]
\includegraphics[height=14cm,width=12cm]{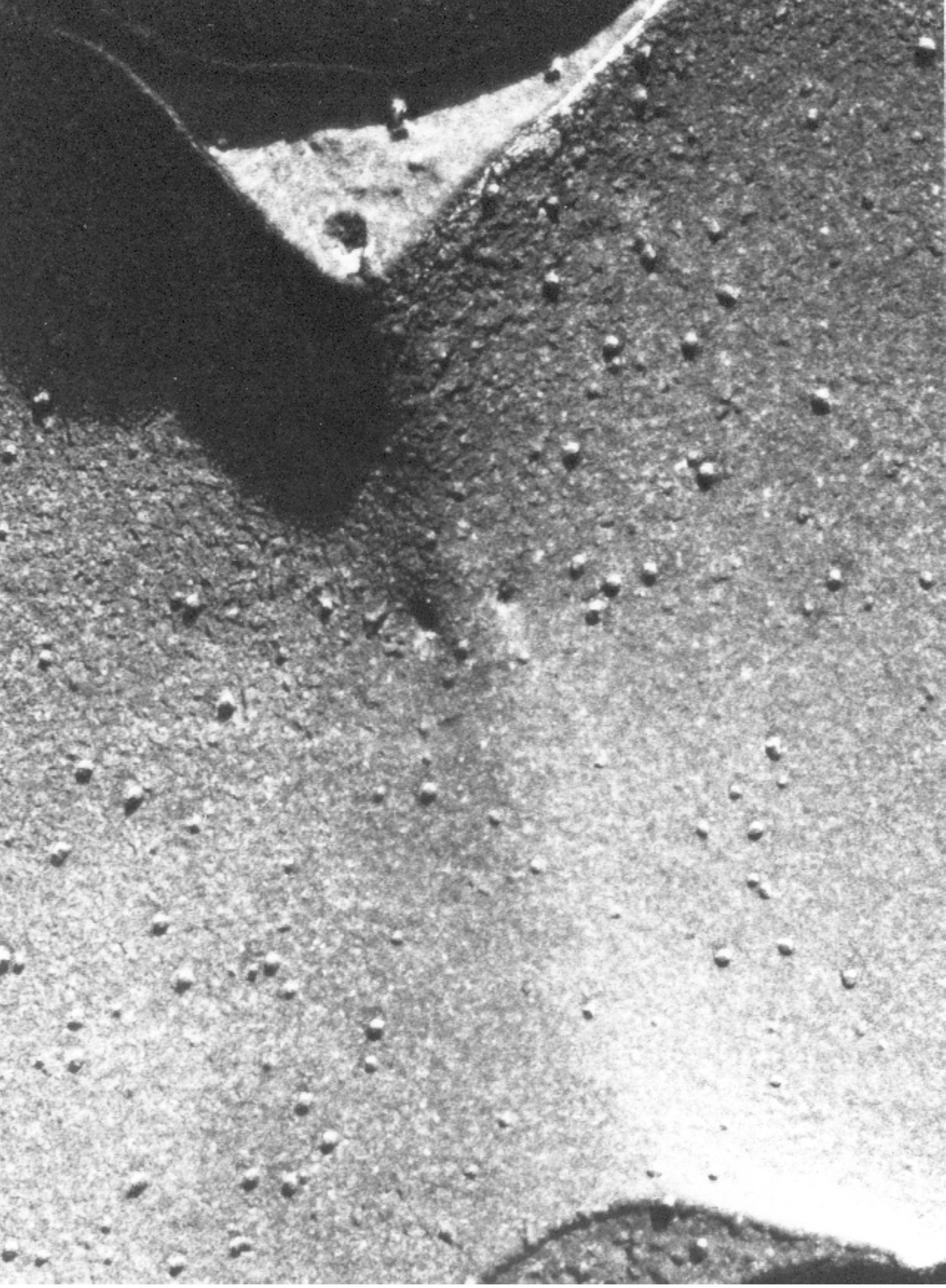}
\caption{Cross section through plasma membrane with dimple and, at
the bottom, a glimpse of the granule. Electron micrograph supplied by
courtesy of D. E. Chandler, Arizona State University, Tempe,
Arizona.} \label{f:pit}
\end{figure}

\subsection{Dimple Formation Prior to the Fusion Event}\label{ss:dimple}

In the introduction to this chapter we defined the fusion pore as the molecular structure
that transiently connects the lumens of two membrane compartments during their fusion.
We emphasized that  making the fusion pore plays a key role in all intracellular trafficking and endocytotic
and exocytotic pathways in all eukaryotic cells, including the regulated exocytosis in endocrine, exocrine and neuronal
cells like our $\beta$-cell. However, from Monck and Fernandez \cite[{\em 1992}]{MonFer} to Shillcock and Lipowski
\cite[{\em 2006}]{ShiLipComp}, researchers agreed that {\em despite its importance, the nature of the fusion pore is unknown}
(\cite[p. 1395]{MonFer}).

\begin{figure}[htb]
\includegraphics[height=10cm,width=12cm]{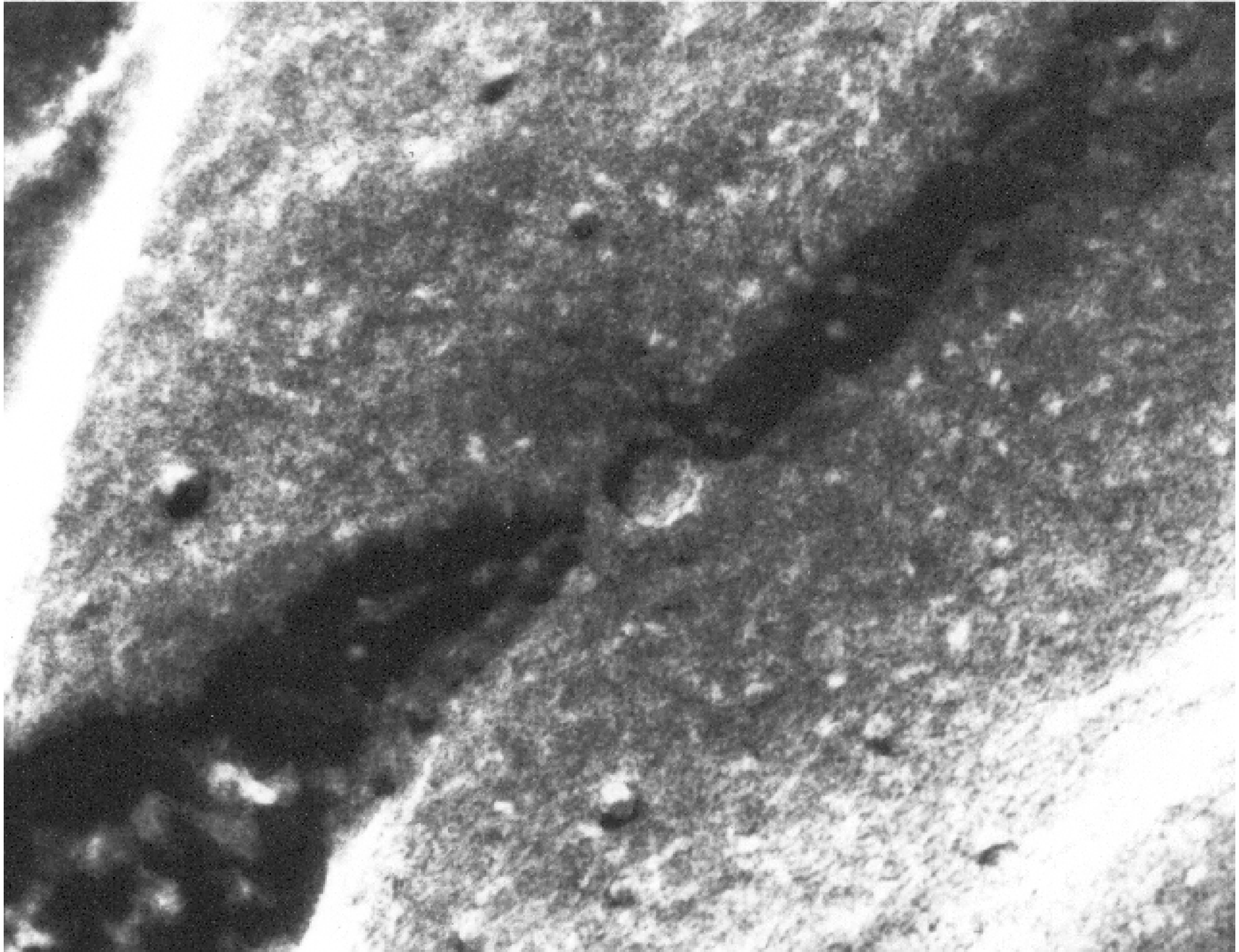}
\caption{The fusion pore shown as imprint of the plasma membrane
dimple onto the granule membrane, from \cite[Fig. 9C]{CCCMZ}.
Reprinted with permission by Springer-Verlag. The high resolution
micrograph was supplied by courtesy of D. E. Chandler, Arizona State
University, Tempe, Arizona.} \label{f:hole}
\end{figure}

In a remarkable series of micrographs, based on rapid freezing techniques for
electron microscopy,
the renown expert on mammal
egg cells, Douglas Chandler from Arizona State University and collaborators captured the formation of the fusion pore in mast cells already some 30 years ago, \cite{{ChaHeu79}, {ChaHeu80}, {CCCMZ}}. They demonstrated
that the pores are made of a curved bilayer which spans the granule and plasma membrane. The micrographs
also gave a hint of the events preceding the making of the fusion pore:  namely the formation of a dimple that
approaches the granule membrane after stimulation, see Figures \ref{f:pit} and \ref{f:hole}.

Since then, the dimple formation has been observed in many different cell types
immediately upon stimulation before
the fusion event, see, e.g., Jena and collaborators \cite{JCJSA}. Using atomic force microscopy (AFM) they
demonstrated the presence of many simultaneous dimples in pancreatic acinar cells after exposing them to a
secretagogue. The paper contains references to analogous demonstrations revealing the presence of pits and depressions
also in pituitary and chromaffin cells prior to secretion.
To the best of our knowledge, capturing the very fusion event of regulated exocytosis in
pancreatic $\beta$-cells has not yet been achieved by present imaging methods.

Not surprisingly, the secretory granules have never been seen to form dimples on their own membrane,
in accordance with the common perception that the plasma membranes are relatively slack and the
membranes of densely packed granule under tension. That corresponds to our preceding Working Hypothesis No. 2.

\subsection{The Flickering of Regulated Exocytosis}\label{ss:flickering}

Another feature of the fusion pore making requires explanation (and
is nicely explained by instabilities of the discussed AC current),
namely the flickering of the fusion event, i.e., the common
observation that the fusion pore can be maintained only after a
while of opening and is not stable immediately after its making.
Irregular rapid pore openings and closures are observed that last
from a few milliseconds to many seconds, see Fernandez et al.
\cite{FNG} for fusion-pore flickering ({\em kiss-and-run}) in mast cells,
and Rosenheck
\cite{Ros} and Jahn et al. \cite{JLS} for wider reports
on the observed flickering of the fusion event, mostly for synaptic vesicle
exocytosis.

\section{The Model}\label{s:model}

In this section we fix our notation and introduce the basic equations for the propagation
of the electromagnetic (wandering) field wave and the making of the dimple.

\subsection{The Force Balance Equation}\label{ss:3-1}

Let $\vec{\textbf{r}}=\vec{\textbf{r}}(x,y,z,t)$ denote the displacement vector of the dimple and $m$ the
mass of an elementary unit of the dimple.
Then the resulting force for making the dimple is approximately equal to $m\dfrac{\partial^2 \vec{\textbf{r}}}{\partial t^2}$, i.e, we can write
\begin{equation}m\frac{\partial^2 \vec{\textbf{r}}}{\partial t^2}=\vec{\textbf{F}}_{\elas}+\vec{\textbf{F}}_{vis}+\vec{\textbf{F}}_{\ext},
\label{Newton-Law}
\end{equation}
where $\vec{\textbf{F}}_{\elas}$ denotes the restoring (or elastic) force, $\vec{\textbf{F}}_{vis}$ stands for viscosity reaction of the medium surrounded our plasma membrane, and $\vec{\textbf{F}}_{\ext}$ is an external force providing the membrane displacement from equlibrium state.

Our first task is to describe the restoring force $\vec{\textbf{F}}_{\elas}$\,.
In first approximation, we shall assume that the plasma membrane is a surface in $\mathbb{R}^3$ without bending resistance. Hence $\vec{\textbf{F}}_{\elas}$ is defined only by the surface tension and the variation of the membrane surface area.

To describe  $\vec{\textbf{F}}_{\elas}$ and $\vec{\textbf{F}}_{\vis}$ more precisely,
we adapt the standard text book model for the suspended vibrating string, respectively,
vibrating plate (see, e.g., Churchill \cite[Sect. 93]{Chu} and Logan \cite{Log}). We
consider $1D$- and $2D$-membranes separately, because the restoring force behaves differently in one-dimensional and many-dimensional cases.  We begin with the $1D$-case, since all the arguments are simpler in that situation.

\subsection{The 1D-Case}\label{ss:3-2}
In the equilibrium state, let our $1D$-membrane coincide with the $x$-axis; let $u=u(x,t)$ denote the displacement  of our plasma membrane from equilibrium state at the point $x$ and at the time $t$; and
let $\rho=\rho(x)$ denote the linear membrane density at the point $x$. We restrict ourselves to sufficiently small deformations, so, for now, we will neglect all the terms that are of higher infinitesimal order with respect to $\dfrac{\partial u}{\partial t}$.

\begin{figure}[htb]
\sidecaption
\includegraphics[scale=.75]{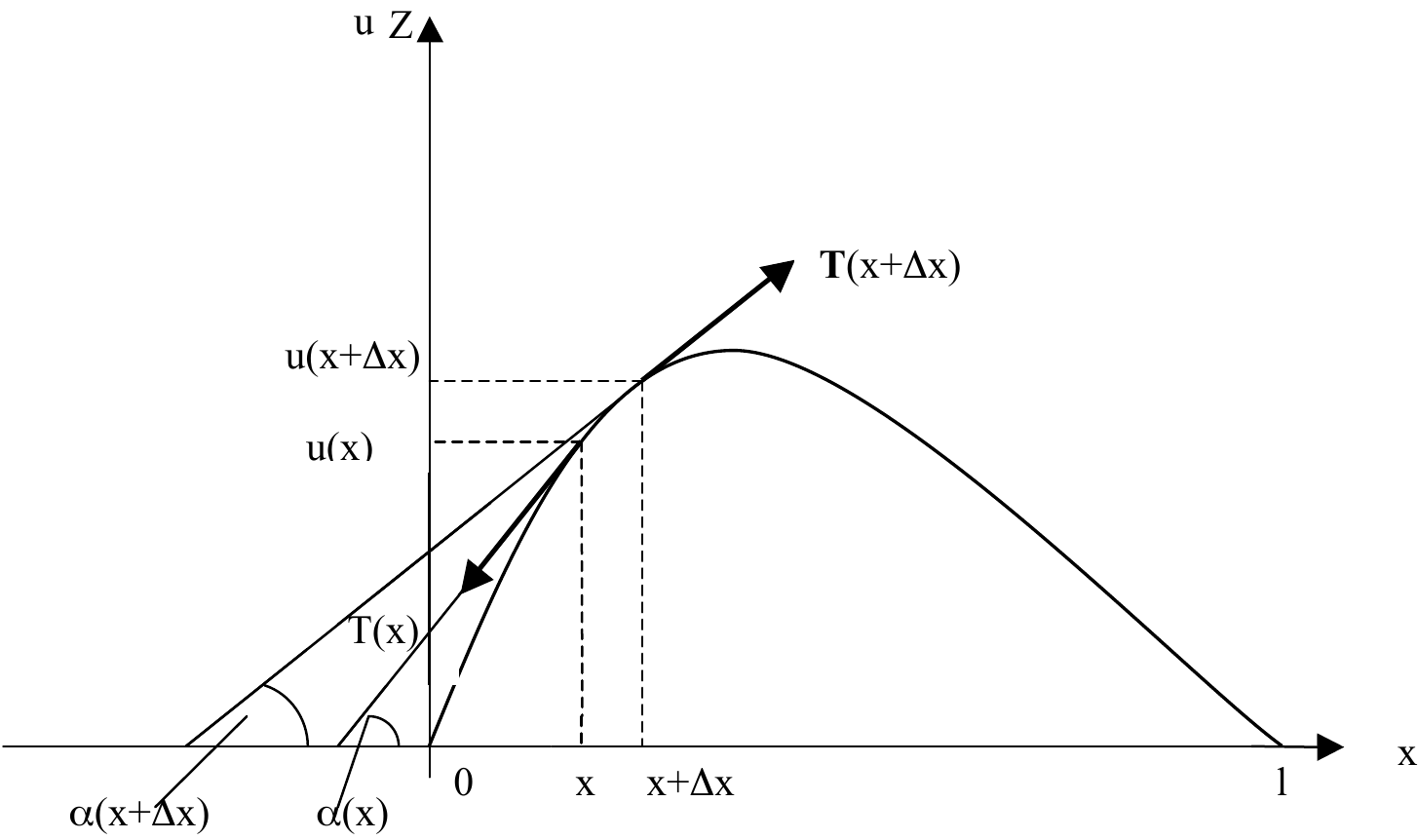}
\caption{Tension along the plasma membrane} \label{f:Fig1}
\end{figure}

Since our plasma membrane has no bending resistance, its tension $\vec{\textbf{T}}(x,t)$ at the point $x$ and at the time $t$ is directed along the tangent to the membrane at $x$. Therefore, the unit $(x, x+\Delta x)$ is subjected to tensions $\vec{\textbf{T}}(x+\Delta x,t)$ and $-\vec{\textbf{T}}(x,t)$ (see Fig. \ref{f:Fig1}). Moreover, according to Hooke's law, $|\vec{\textbf{T}}(x,t)|$ does not depend on $x$ and $t$, i.e., $|\vec{\textbf{T}}(x,t)|=\textbf{T}_0$.

Defining $T_1$ as a relaxation time due to the action of the surrounded medium, we note that
$\vec{\textbf{F}}_{\vis}$ is directed parallel to the vertical coordinate axis and can be modeled as being proportional to $\dfrac{1}{T_1}\dfrac{\partial u}{\partial t}$.

Let us denote by $f(x,t)$ the density of an external force $\mathbf{F}_{\ext}$ acting on the membrane point $x$ at the time $t$ and directing along the vertical axis.
For a description of our electromagnetic candidate, see below Sect. \ref{ss:lorentz}.

Thus, projecting Eq. \eqref{Newton-Law} onto the vertical coordinate
axis and taking the preceding relations into account, we get the
equality

\begin{multline} \label{1a}
\rho \Delta x \frac{\partial^2 u}{\partial t^2}\\=T(x+\Delta x,t)
\sin{\left( \alpha (x+\Delta x)\right) }-T(x,t)\sin{\alpha (x)}
-\frac{2 \rho}{T_1} \Delta x \frac{\partial u}{\partial t}+f(x,t)\Delta x.
\end{multline}
In the context of our approximation
$$
\sin{(\alpha )}=\frac{\tan{(\alpha )}}{\sqrt{1+\tan^2{(\alpha )}}}\approx \tan{(\alpha ) }=\frac{\partial u}{\partial x},
$$
and, consequently, equality (\ref{1a}) takes the form
\begin{equation}
\frac{\partial^2 u}{\partial t^2}=\frac{\textbf{T}_0}{\rho \Delta x}\left[ \frac{\partial u(x+\Delta x,t)}{\partial x}-
\frac{\partial u(x,t)}{\partial x}\right] -\frac{2}{T_1}\frac{ \partial u}{\partial t} +\frac{1}{\rho}f(x,t).
\label{1b}
\end{equation}
Passing in (\ref{1b}) to the limit as $\Delta x \to 0$ we arrive at
\begin{equation}
\frac{\partial^2 u}{\partial t^2}+\frac{2}{T_1}\frac{ \partial u}{\partial t} -c_s^2
\frac{\partial^2 u}{\partial x^2} = \frac{1}{\rho}f(x,t),
\label{2a}
\end{equation}
where $c_s^2=\dfrac{\textbf{T}_0}{\rho}$ is the speed of sound in the membrane.

\subsection{The 2D-Case}\label{ss:3-3}

\begin{figure}[htb]
\sidecaption
\includegraphics[scale=.85]{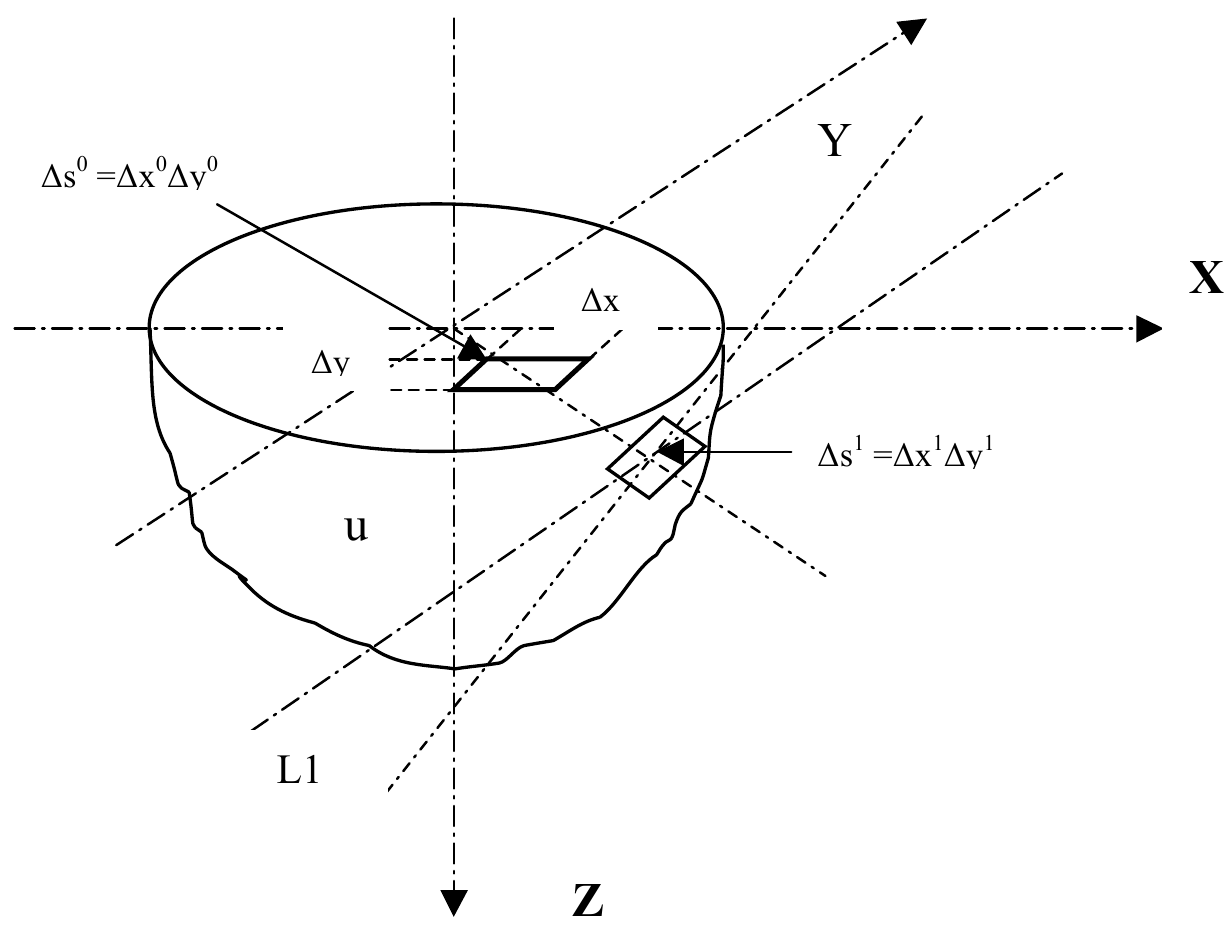}
\caption{Notations for membrane displacement from equilibrium state at the point
$(x,y)$ and at the time $t$} \label{f:Fig2}
\end{figure}

Similarly to the $1D$-case we assume that in equilibrium state our plasma membrane lies in a subspace $XY$
(see Fig. \ref{f:Fig2}) and $u=u(x,y,t)$ denotes the membrane displacement from equilibrium state at the point
$(x,y)$ and at the time $t$. We will consider only small deformations such that
\begin{equation}
\left( \frac{\partial u}{\partial x}\right)^2 \ll 1, \quad  \left( \frac{\partial u}{\partial y}\right)^2 \ll 1.
\label{small-deformations}
\end{equation}

Let $d \sigma$ be the unit length of some closed path lying on the membrane surface, and let $P$ be a point belonging to $d\sigma$. Then the unit $d\sigma$ is subjected to the tension force $\vec{\textbf{T}}d\sigma$, where $\vec{\textbf{T}}=\vec{\textbf{T}}(x,y,t)$ denotes the surface tension. Due to absence of the membrane resistance to bending and shear, we can say that the vector $\vec{\textbf{T}}$ always lies on the hyperplane $L_1$ tangential
to the membrane surface at the point $P$, and $\vec{\textbf{T}}$ is orthogonal to $d\sigma $ (see Fig.~\ref{f:Fig2}).

 In addition, inequalities (\ref{small-deformations}) guarantee that the tangential hyperplane $L_1$ lies almost parallel to the hyperplane $XY$. To prove this statement, we have to show that the length of the
projection of the vector $\vec{\textbf{T}}(x,y,t)$ onto $XY$ is approximately equal
to $|\vec{\textbf{T}}(x,y,t)|$. Indeed, by definition
$$
\bigl|\text{Projection}_{XY}\vec{\textbf{T}}(x,y,t)\bigr|=|\vec{\textbf{T}}(x,y,t)|\cos{(\beta)},
$$
 where $\beta$ stands for the angle between the tension vector $\vec{\textbf{T}}$ and hyperplane $XY$.
It is easy to see that $\beta$ is not bigger than the angle $\gamma$ between the tangent hyperplane $L_1$ and $XY$. Therefore,
$$
\cos{(\beta)} \geqslant \cos{(\gamma)}=\left[ 1+\left( \frac{\partial u}{\partial x}\right)^2+\left(\frac{\partial u}{\partial y}\right)^2\right]^{-1/2}\approx 1,
$$
and, consequently,
$$
\bigl|\text{Projection}_{XY}\vec{\textbf{T}}(x,y,t) \bigr|\approx |\vec{\textbf{T}}(x,y,t)|.
$$
Hooke's law guarantees that $|\vec{\textbf{T}}(x,y,t)|$ does not
depend on the $t$-variable, whereas the orthogonality of
$\vec{\textbf{T}}(x,y,t)$ and $d \sigma$ provides the independence
of $|\vec{\textbf{T}}|$ on the variables $x$ and $y$ as well. It
means that
$$
|\vec{\textbf{T}}(x,y,t)|=\textbf{T}_0=\text{const}.
$$

Now, considering a rectangle unit $ds=\Delta x \Delta y$ on the membrane surface, we can write the restoring force acting on this unit as:
\begin{multline*}
 \textbf{T}_0 \Delta y \left[\left( \frac{\partial u}{\partial x}\right) \bigg|_{x+\frac{\Delta x}{2}}  -\left( \frac{\partial u}{\partial x}\right) \bigg|_{x-\frac{\Delta x}{2}} \right] +
 \textbf{T}_0 \Delta x \left[\left( \frac{\partial u}{\partial y}\right) \bigg|_{y+\frac{\Delta y}{2}}  -\left( \frac{\partial u}{\partial y}\right) \bigg|_{y-\frac{\Delta y}{2}} \right]\\
=\textbf{T}_0 \Delta y \frac{\partial^2 u}{\partial x^2}\Delta x+\textbf{T}_0 \Delta x \frac{\partial^2 u}{\partial y^2}\Delta y
= \textbf{T}_0 \left( \frac{\partial^2 u}{\partial x^2}+\frac{\partial^2 u}{\partial y^2}\right) \Delta x\Delta y.
\end{multline*}

It remains to describe the external and viscosity forces acting on $ds$.
Similarly to the $1D$-case, $f(x,y,t)$ denotes the density of the external force $\mathbf{F}_{\ext}$ at the point $x$ and at the time $t$. It
 is directed orthogonally to the membrane surface, while $\vec{\textbf{F}}_{\vis}$, directed opposite to the vector $\vec{\textbf{F}}_{\ext}$, is proportional to the velocity $\dfrac{\partial u}{\partial t}$.

Let $\rho (x,y)$ denote the membrane surface density, then the mass of the unit $ds$ is equal $\rho (x,y)\Delta x\Delta y$. Finally, defining  a relaxation time due to the action of the surrounded medium by $T_1$, we can write the variant of Eq. (\ref{Newton-Law}) for $2D$-membranes as follows:
\begin{equation}\rho \Delta x\Delta y \frac{\partial^2 u}{\partial t^2}=\textbf{T}_0\left( \frac{\partial^2 u}{\partial x^2}+\frac{\partial^2 u}{\partial y^2}\right) -\frac{2\rho}{T_1}\Delta x\Delta y\frac{\partial u}{\partial t}+f(x,y,t)\Delta x \Delta y.
\label{3a}
\end{equation}
After elementary transformations, Eq. (\ref{3a}) takes the form
\begin{equation}
\frac{\partial^2 u}{\partial t^2}+\frac{2}{T_1}\frac{ \partial u}{\partial t} -c_s^2 \left(  \frac{\partial^2 u}{\partial x^2}
+\frac{\partial^2 u}{\partial y^2}\right)
=\frac{1}{\rho}f(x,y,t), \qquad c_s^2=\frac{\textbf{T}_0}{\rho}.
\label{3}
\end{equation}

\subsection{Further Approximations}\label{ss:approximations}

Having equations (\ref{2a}) and (\ref{3}) at hand, we observe that the process of dimple forming is quasi-static, i.e.,
$$
\frac{\partial^2 u}{\partial t^2}\ll 1.
$$
The latter means that we can neglect this term in both equations.

It should also be pointed out that $u=|Du|=0$ at those membrane points where there is no influence of the  external force. Here $Du$ denotes the spatial gradient of the displacement $u$.

Recall that the characteristic function $\chi_E$ of a set $E$ is defined by
$$
\chi_E(z)=\left\lbrace
\begin{aligned}
&1, \quad \text{for}\ z\in E,\\
&0, \quad \text{for}\ z \notin E.
\end{aligned} \right.
$$

Taking into account all the above remarks we get the following model equations for the dimple forming
in the one- and two-dimensional cases, respectively:
\begin{align}
\frac{2}{T_1}\frac{\partial u}{\partial t}-c_s^2 \frac{\partial^2 u}{\partial x^2}&
=\frac{1}{\rho}f(x,t)\chi_{\{u>0\}},
\label{1d-model}\\
\frac{2}{T_1}\frac{\partial u}{\partial t}-c_s^2 \left( \frac{\partial^2 u}{\partial x^2}+\frac{\partial^2 u}{\partial y^2}\right)&=
\frac{1}{\rho}f(x,y,t)\chi_{\{u>0\}}.
\label{2d-model}
\end{align}


\subsection{Lorentz Force}\label{ss:lorentz}

Clearly, there is a variety of external forces resulting from the
electromagnetic field wave. One can expect that they all play
together in forming the dimple. However, taking our point of
departure in an alternating current, we shall concentrate on the
{\em Lorentz force}, (implicitly deal with the {\em Coulomb force}), but discard the {\em dipole electric
force}, the {\em magnetic force}, and the {\em van der Waals force} for now.

\subsubsection{Peculiarity of the Lorentz Force}
In our physiological context, the peculiar role of the Lorentz force
is that it exerts its action in one fixed direction, the direction of the
propagation of the field wave, even when it is related to an alternating
current. Roughly speaking, it is the sophistication and power of
the electromagnetic aspects of the regulated exocytosis that a
relatively weak and extremely low frequent electrodynamic wave can
transport energy along a straight line over a large intracellular
distance and exert its action on the charged phospholipid molecules in
the plasma membrane. These charged molecules make a substantial
part of the eukaryotic plasma membrane, around 11\% according to Rosenheck \cite{Ros}.

 \begin{figure}[htb]
\includegraphics[scale=.68]{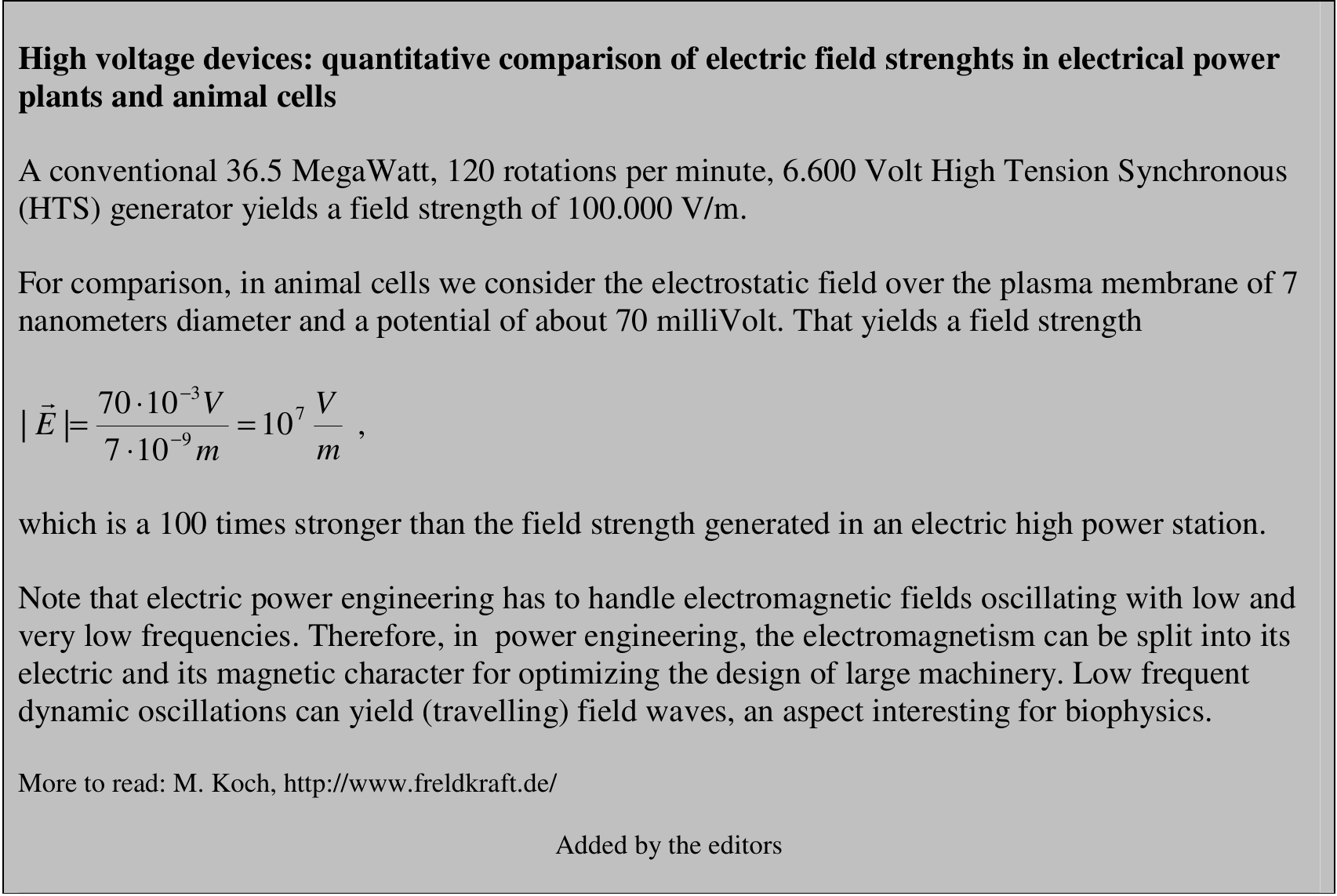}
\end{figure}

\subsubsection{Energy Estimates} It may be instructive to have a rough idea of the scale
of the electric forces and related energies around the dimple
formation. see also the box. Based on the estimates given by
Rosenheck l.c., we have a surplus of about $10^{10}$ (negative)
charges per 1000 $\nm^2$ plasma membrane area, i.e., about $10^9$
charge carriers around the dimple top stretching over an area of
approximately 100 $\nm^2$\,. The energy of one charged phospholipid
molecule was calculated by Rosenheck l.c. as being $10^{-19}
\operatorname{W}\operatorname{sec} = 10^{-19}
\operatorname{N}\operatorname{m}$. Dividing by a characteristic
length of $10 \nm$ we obtain a force of around $10^{-11}
\operatorname{N}$, i.e., a total force of $10^{-2} \operatorname{N}$
exerted on the dimple region.

For comparison, the gravitational force of the mass $10^{-21} \operatorname{kg}$
of a dimple of about $10^3\nm$ volume
(and specific weight comparable to water) would only sum up to $10^{-21}\times 9.8
\sim 10^{-20}\operatorname{N}$.
For another comparison, we refer to the Koch-Stetter electromagnetic field wave generator. For now,
its B-field is of about $35 \operatorname{mT}$ (milliTesla) and the exerted forces are
visible and can be measured.

To sum up, for the making of the dimple we are investigating electric forces which surpass by large the above estimated $10^{-2}\operatorname{N}$.

\subsubsection{Model Quantities} Our first model quantity is the capacitive
reactance $X_C$ as scalar function of place and time, depending on the
shape of the dimple. Our second model quantity
is the time dependent (namely dimple shape and $X_C$ dependent) amplitude $\widehat
{\mathbf B}$ of the magnetic flux density.

\subsubsection{Model Equations} As mentioned before, we have two model levels.

(I) There must be a variational equation,
minimizing an energy functional or another related cost functional
which gives the change and spatial distribution of $X_C$.

(II) As discussed above, there are also elastic and viscous forces resisting
the necessary re-packing of the lipid heads under deformation
and pulling the plasma membrane back
in the more smooth non-dimple form, interrupting exocytosis when the
electrical and magnetic production is interrupted.

\subsubsection{The Lorentz Force} Our favorite external force for
making the dimple is the Lorentz force
\begin{equation}\label{e:lor}
{\mathbf F}_L =  q\mathbf{E}_0 + q \bigl(\mathbf{v}\times
\mathbf{B}\bigr) - \gamma \mathbf{v},
\end{equation}
where $q$ is the charge (i.e., twice the number) of the released
$\Ca^{2+}$ ions, ${\mathbf E}_0$ is the background electric field,
here assumed to be zero, $\mathbf B$ is the magnetic flux density,
$\mathbf v$ is the velocity of the charged dimple particles with
respect to the coordinate system in which the quantities
$\mathbf{F},\ \mathbf{E}_0$, and $\mathbf{B}$ are calculated, i.e.,
$\mathbf {v}$ is up to the sign the wandering velocity of the
directed electromagnetic field wave, and $\gamma\mathbf{v}$ is the
friction force under propagation. Once again, equation \eqref{e:lor} is written
in SI. In the cgs system, common in physics, the term $q
\bigl(\mathbf{v}\times \mathbf{B}\bigr)$ has to be divided by $c$,
the velocity of light in vacuum.

Note that we have for a (moving) field wave (see, e.g., \cite[Chapter 13]{Ti})
\begin{align}\label{e:b}
\mathbf{B}(\mathbf{x},t) &=\widehat {\mathbf
B}\cos\bigl(\mathbf{k}\cdot\mathbf{x}(t)- k\,
v(\mathbf{x},t)\, t\bigr)
 +\mathbf{B}_{\operatorname{DC}}\,,\\
\label{e:e} \mathbf{E}(\mathbf{x},t)
&=\mathbf{v}(\mathbf{x},t)\times\mathbf{B}(\mathbf{x},t)\,.
\end{align}
Here, $\mathbf{B}_{\operatorname{DC}}$ denotes the background B-field corresponding
to the direct current $\mathbf E_0$\,, $\mathbf x = (x,y,z)$ denotes a position,
$\mathbf k$ denotes the wave vector with $k=|\mathbf{k}|$
and $v(\mathbf{x},t)= |\mathbf{v}(\mathbf{x},t)|$.

The Lorentz force of Eq. \eqref{e:lor} can be inserted into the general
balance equation \eqref{Newton-Law}, into the 1-dimensional
model equations \eqref{1a}, \eqref{2a}, simplified to \eqref{1d-model},
and into the
2-dimensional model equations \eqref{3a}, simplified to \eqref{2d-model}.

\subsubsection{Work Equation} In Sect. \ref{ss:magnetic-closure}, we have sketched the feedback mechanism  of
the dimple forming and explained why and how the forming of the dimple strengthens the
electromagnetic field wave. It is beyond the range of this chapter to model
that mechanism in detail as a free boundary problem. Basically, to relate our defining quantity $X_C$ to the listed
balance equations, we shall express the power in terms of $X_C$ and then derive an integral for the
electromagnetic energy density (per volume), where $X_C$ enters. The details will be worked out
separately.

\section{Apposite Results on Parabolic Obstacle Problems}\label{s:regularity}

At the end of Sect. \ref{ss:approximations}, we introduced the model equations
\eqref{1d-model} and \eqref{2d-model} that can be treated as parabolic free boundary problems (FBP).

\subsection{Review of Free Boundary Problems}\label{ss:free-boundary}

The expression \textbf{free boundary problem} means that we deal with a {problem with two a-priori unknown objects}: an unknown {\em set} coming up in a {\em solution} of a partial differential equation.

A typical example is the {\em Stefan problem} describing the melting
of an ice cube in a glass of water. If ice begins to melt then the
region occupied by water will grow and the interface-surface between
the ice and the water (it is called the free boundary) will move and
change its shape, see, e.g., Friedman \cite[Sect. 1.9]{Fri}. As
another typical example of FBP we mention the flame propagation
problem describing the evolution of the flame front, for that see
the various contributions in \cite[Chapter 8]{FasPri}.

Using a common transformation of the independent variables, we may
normalize the coefficients and so reduce our model equations
\eqref{1d-model} and \eqref{2d-model} to the following problem:
\begin{equation}
\left\lbrace \begin{aligned}
&\Delta u(x,t)-\frac{\partial u}{\partial t}(x,t)=f(x,t)\chi_{\left\lbrace u>0\right\rbrace },\\
&u(x,t)\geqslant 0
\end{aligned} \right. \qquad \text{almost everywhere (a.e.) in } \mathcal{D},
\label{statement-obstacle-problem}
\end{equation}
where $(x,t)$ denotes the points in $\mathbb{R}^n \times \mathbb{R}$ with the space variable $x=(x_1,\dots,x_n)$ belonging to $\mathbb{R}^n$ and the time variable $t$ belonging to $\mathbb{R}$, $\Delta$ is the Laplace operator defined as
$$
\Delta u=\sum\limits_{i=1}^n \frac{\partial^2 u}{\partial x_i^2},
$$
$\chi_{\left\lbrace u>0\right\rbrace }$ is the characteristic
function of the set $\left\lbrace(x,t) \in \mathbb{R}^{n+1}\mid
u>0\right\rbrace $ (see Section 3.4 for the precise definition),
$\mathcal{D}$ is a given open set in $\mathbb{R}^{n+1}$\,, and $u$ is a
locally bounded weak (i.e., in the distributional sense) solution.

Observe that $\left\lbrace u>0\right\rbrace$ is a priori an unknown
open subset of $\mathcal{D}$. We denote by $\Gamma$ the intersection
of $\mathcal{D}$ with the \textit{boundary} of the set $\left\lbrace
u= 0\right\rbrace $. We will call $\Gamma$  the \textit{free
boundary}. For our dimple forming, we have $\Gamma =
\{\Gamma(t),\{t\}\}_{t\ge 0}$ of Fig. \ref{f:vesicle}.

If, additionally, the condition
\begin{equation}
\frac{\partial u}{\partial t} \geqslant 0 \qquad \text{a.e. in}\quad \mathcal{D},
\label{sign-of-ut}
\end{equation}
is satisfied, then our FBP (\ref{statement-obstacle-problem})
becomes the Stefan problem mentioned above. Inequality
(\ref{sign-of-ut}) means that our dimple is formed without
``returning back". In general, however, we can not guarantee
(\ref{sign-of-ut}) from the assumptions given in Sections
\ref{ss:3-2} and \ref{ss:3-3} only.

For the function $f(x,t)$, we assume that
\begin{itemize}
\item[(1)] \quad $f$ is {\em non-degenerate} in $\mathcal{D}$, i.e., there exists $\delta_0>0$ such that $f(x,t)>\delta_0$ for any $(x,t) \in \mathcal{D}$;
\item[(2)] \quad $f$ is {\em H{\"o}lder continuous} in $\mathcal{D}$ with some $\alpha \in (0,1)$, i.e., $f$ is a bounded continuous function in $\mathcal{D}$, and for all points $(x,t)$ and $(y,s) \in \mathcal{D}$ such that $(x,t) \neq (y,s)$ we have the inequality
$$
\frac{|f(x,t)-f(y,s)|}{\left( |x-y|^2+|t-s|^2\right)^{\alpha /2}
}<\infty.
$$
\end{itemize}

\subsection{Qualitative Properties of Solutions}

Let $u$ be a solution of FBP (\ref{statement-obstacle-problem}),
let $z^0=(x^0,t^0)\in \mathcal{D}$, let $\rho>0$ sufficiently small,
and let $Q_{\rho}(z^0):=\left\lbrace |x-x^0|<\rho\right\rbrace \times (t^0-\rho^2, t^0)$.
The following estimates provide us with information about the behaviour of our displacement $u$ near the interface between the sets $\left\lbrace u>0\right\rbrace $ and $\left\lbrace u=0\right\rbrace$, i.e., near the free boundary. 

\begin{itemize}
\item \quad There exists a constant $c=c(n)>0$ such that
\begin{equation}\label{e:estimate1}
\sup\limits_{Q_{\rho}(z^0)} u \geqslant c(n) \rho^2.
\end{equation}
This nondegeneracy estimate holds true  for all points $z^0$
belonging to the closure of the set $\left\lbrace u>0\right\rbrace $
and for all $\rho$ sufficiently small.
\end{itemize}

Moreover, our solution $u$ has quadratic growth near the free
boundary:
\begin{itemize}
\item \quad There exists a constant $C>0$ such that
\begin{equation}\label{e:estimate2}
\sup\limits_{Q_{\rho}(z^0)} u \leqslant C\rho^2.
\end{equation}
This nondegeneracy estimate holds true  for all free boundary points $z^0\in \Gamma$ and for all $\rho$ sufficiently small.
\end{itemize} 

\begin{itemize}
\item \quad There exists a universal constant $M>0$ such that
\begin{equation}\label{e:estimate3}
\sup\limits_{Q_{\rho}(z^0) \cap \{u>0\}}\left(
\left|\frac{\partial^2 u}{\partial x_i \partial x_j}\right|+
\left|\frac{\partial u}{\partial t}\right|\right) \leqslant M .
\end{equation}
This inequality holds true for all free boundary points $z^0 \in \Gamma$  and for all $\rho$ sufficiently small.
\end{itemize} 

For $f(x,t)=const$ all these three statements were proved in
\cite{CPS}. The case of general $f=f(x,t)$ was considered in \cite{Bl2}
for $n>1$ and in \cite{BDM} for $n=1$. 

The preceding estimates \eqref{e:estimate1}, \eqref{e:estimate2},
and \eqref{e:estimate3} indicate that characteristic base diameter,
growth rate, depth and time of the dimple forming are mathematically
well defined and, therefore, these values should, in principle, be
measurable in experiment, see below Sect. \ref{ss:suggested}.

\subsection{Classification of Blow-up Limits in $\mathbb{R}^{n+1}$}

The idea is to use \textit{blow-up} sequences, which are a kind of zooms, and to
look at the "infinite zoom".  Suppose that $u$ is a solution of
the problem (\ref{statement-obstacle-problem}), $z^*=(x^*,t^*)$ is a
free boundary point and $f(x^*,t^*)\neq 0$. For $\lambda >0$\, consider the functions
$$
u_{\lambda}(x,t):= \frac{u\left(x^*+x\frac{\lambda}{\sqrt{f(x^*,t^*)}}, t^*+t\frac{\lambda^2}{f(x^*,t^*)}\right)}{\lambda^2}, \quad \text{for} \quad (x,t) \in \mathcal{D}_{\lambda}:=\frac{1}{\lambda}\mathcal{D}.
$$
There exists a sub-sequence $\left\lbrace \lambda_k\right\rbrace $
converging to zero such that the blow-up sequence $\left\lbrace
u_{\lambda_k}\right\rbrace $ converges to one of the following
blow-up limits:

\begin{itemize}
\item \quad $u_0=u_{0,e}(x,t):=\frac{1}{2}(x^T \cdot \vec{\bf{e}})^2$, for a unit vector $\vec{\bf{e}}$, where $x^T \cdot \vec{\bf{e}}$ denotes the scalar product in $\mathbb{R}^n$,
\item \quad $u_0=u_{0,m}(x,t)=mt+x^T\cdot \mathcal{M}\cdot x$, where $m$ is a constant and $\mathcal{M}$ is a $(n\times n)$-matrix satisfying
$\textbf{Tr}\mathcal{M}=m+1$.
\end{itemize}

Observe that the blow-up limits can (in general) depend on the
choice of the sub-sequence $\left\lbrace \lambda_k\right\rbrace $.
But  it should be emphasized that in view of the non-negativity of
$u$, the limit function $u_0$ is the unique non-negative
distributional solution of
$$
\Delta u -\frac{\partial u}{\partial t}=\chi_{u>0} \qquad \text{a.e. in} \quad \mathbb{R}^n\times (-\infty, t^*).
$$
This means that in the second case $m$ and $\mathcal{M}$ are defined uniquely as well. \vspace{0.3cm}

For $f(x,t)=const$ these results were obtained in \cite{CPS}. For the general case $f=f(x,t)$ we refer the reader to $\cite{Bl2}$.

\subsection{Classification of the Free Boundary
Points}\label{ss:class_gamma}

Going to the regularity properties of the free boundary we observe that

\begin{itemize}
\item \quad The free boundary is a closed set of zero $(n+1)$-Lebesgue measure.
\item \quad The free boundary $\Gamma=\left\lbrace \text{"regular points"} \right\rbrace \cup
\left\lbrace \text{"singular points"} \right\rbrace$.
\item \quad The set of the singular free boundary points is closed.
\end{itemize}
The {\em singular points} are defined as the free boundary points
for which there exists a blow-up limit of the second type, i.e.,
$u^0(x,t)=mt+x^T\cdot \mathcal{M}\cdot x$. The set $\Gamma\setminus
\left\lbrace \text{"singular points"}\right\rbrace $ is the set of
regular points. For singular points and for
$k=0,\dots,n$ the sets $S(k)$ are considered additionally, where $S(k)$ is
defined  as the set of singular points such that $\text{dim}\,
\text{Kern}\, \mathcal{M}=k$ and the smallest of the $k$ non-zero
eigenvalues is bounded from below by a fixed positive constant.

The complete classification of all free boundary points can be given via a relatively new approach introduced by G.S.\,Weiss in \cite{W}. For a solution $u$ of FBP (\ref{statement-obstacle-problem}) and for a free boundary point $z^*=(x^*,t^*)$ consider the following  energy functional
$$
W(\tau, z^*, u,f):=\frac{1}{\tau^4}\int\limits_{t^*-4\tau^2}^{t^*-\tau^2}\int\limits_{|x-x^*|<\tau}
\left( |\nabla u|^2+2fu+\frac{u^2}{t-t^*}\right) G(x-x^*,t^*-t)dxdt,
$$
where $G(x,t)$ is the heat kernel. 

The functional $W$ has the following remarkable properties:
\begin{itemize}
\item \quad $\lim\limits_{\tau \to 0^+}W(\tau, z^*,u,f)$ exists and is finite;
\item \quad There are only two possible values for $\lim\limits_{\tau \to 0^+}W(\tau, z^*,u,f)$, namely
$$
\lim\limits_{\tau \to 0^+}W(\tau, z^*,u,f)=\left\lbrace
\begin{aligned}
&A_n, \ \quad \text{if}\quad z^* \quad \text{ is a regular point}\\
&2A_n, \quad \text{if}\quad z^* \quad \text{ is a singular point}
\end{aligned}
\right.
$$
for some constant $A_n>0$ depending only on the dimension $n$.
\end{itemize} 

Finally, we observe that
\begin{itemize}
\item Around regular points the free boundary is a smooth graph.
\item Singular points belonging to $S(n)$ are isolated.
\item $S(n)$ is contained locally in a $C^2_x$-graph in space.
\item $S(k)$ for $0 \leqslant k \leqslant n-1$ is contained locally in a $k$-manifold of class $C^{1/2}_{x,t}$.
\end{itemize} 

For all results concerning the regular points we refer to \cite{CPS}. The results about singular points were proved in \cite{Bl2} (see also \cite{Bl1}).

\section{Conclusions}\label{s:drive}

Our conclusions consist of preliminary findings which will require further experiments and measurements
to be confirmed - or falsified.

\subsection{Summary of (Partly Speculative) Working Hypotheses}
For better reading, we begin this closing section by summarizing our assumptions and choices of
emphasis.

\begin{enumerate}
\item Among all relevant aspects of the well- and malfunctioning of pancreatic $\beta$-cells,
we focus on a single membrane process, the lipid bilayer fusion event.
\item We suppose that future imaging will prove the making of a plasma membrane dimple
before the making of the fusion pore.
\item We assume that the dimple making is essential for the performance of regulated exocytosis
also in $\beta$-cells.
\item We suspect that maintaining the fusion pore and continuing release of the content of the
insulin granules will be interrupted when the dimple is not preserved.
\item We argue for a long-distance regulation of the dimple making (and, hence,
the secretion process), i.e., we claim that the docking of readily
releasable insulin granules at the plasma membrane and the
consecutive making of the fusion pore and the release of the hormone
molecules is induced not as a purely local phenomenon and hence can
not be explained solely by glucose stimuli and corresponding
$\Ca^{2+}$ influx through ion channels in the neighborhood of the
release site. Instead of that, we claim that stimuli and ion influx
also far from the exocytosis site and mediated by intra-cellular
signaling and energy transport will be decisive for initiating and
maintaining regulated exocytosis even at one isolated site.
\item We point to a low-frequent electromagnetic field wave as a possible regulator of exocytosis.
\item We suppose that the electromagnetic field will be closed over the plasma membrane due to
the iron content of enzymes embedded between the phospholipids of the plasma membrane.
\item We assume that the electrical activity of the mitochondria is equally important as the electrical activity at the plasma membrane. We suppose that the mitochondria
sequester and release $\Ca^{2+}$ ions in a self-regulated way which builds an electromagnetic field and generates a directed (travelling) field wave.
\item We suppose that the synchronization of the electrical activity of
neighboring mitochondria is due to energy efficiency.
\item While many aspects of $\beta$-cell function will require an analysis of their
collective functioning in the Langerhans islets, we conjecture that essential aspects of
regulated exocytosis can be observed on the level of a single cell.
\end{enumerate}

\subsection{The Findings}\label{ss:findings}

We have provided a mathematical model for the initiation of
regulated exocytosis and the making of the fusion pore. The model
relates the geometry and the dynamics of one single membrane
process, namely the forming of an inward oriented dimple in the
plasma membrane before the fusion event, with electromagnetic
features of intracellular calcium oscillations. The model suggests a
new explanation for the observed flickering of regulated exocytosis,
the vanishing of the first phase of secretion in stressed or tired
$\beta$-cells and the final halt of all secretion in overworked
dysfunctional cells: the electromagnetic free boundary model points to
the lack of stability and coordination of the intracellular
$\Ca^{2+}$ oscillations prior to the bilayer membrane vesicle
fusion. We recall that the field character of these oscillations is
magnetic (therefore transferring energy to the fusion site at the
plasma membrane without any  loss). It must be distinguished from
the widely studied $\Ca^{2+}$ influx changing the electrostatic
potential across the plasma membrane and accompanying regulated
exocytosis.

The model is based solely on physical First Principles. All parameters have
a biophysical meaning and can, in principle, be measured.

 \begin{figure}[htb]
\includegraphics[scale=.68]{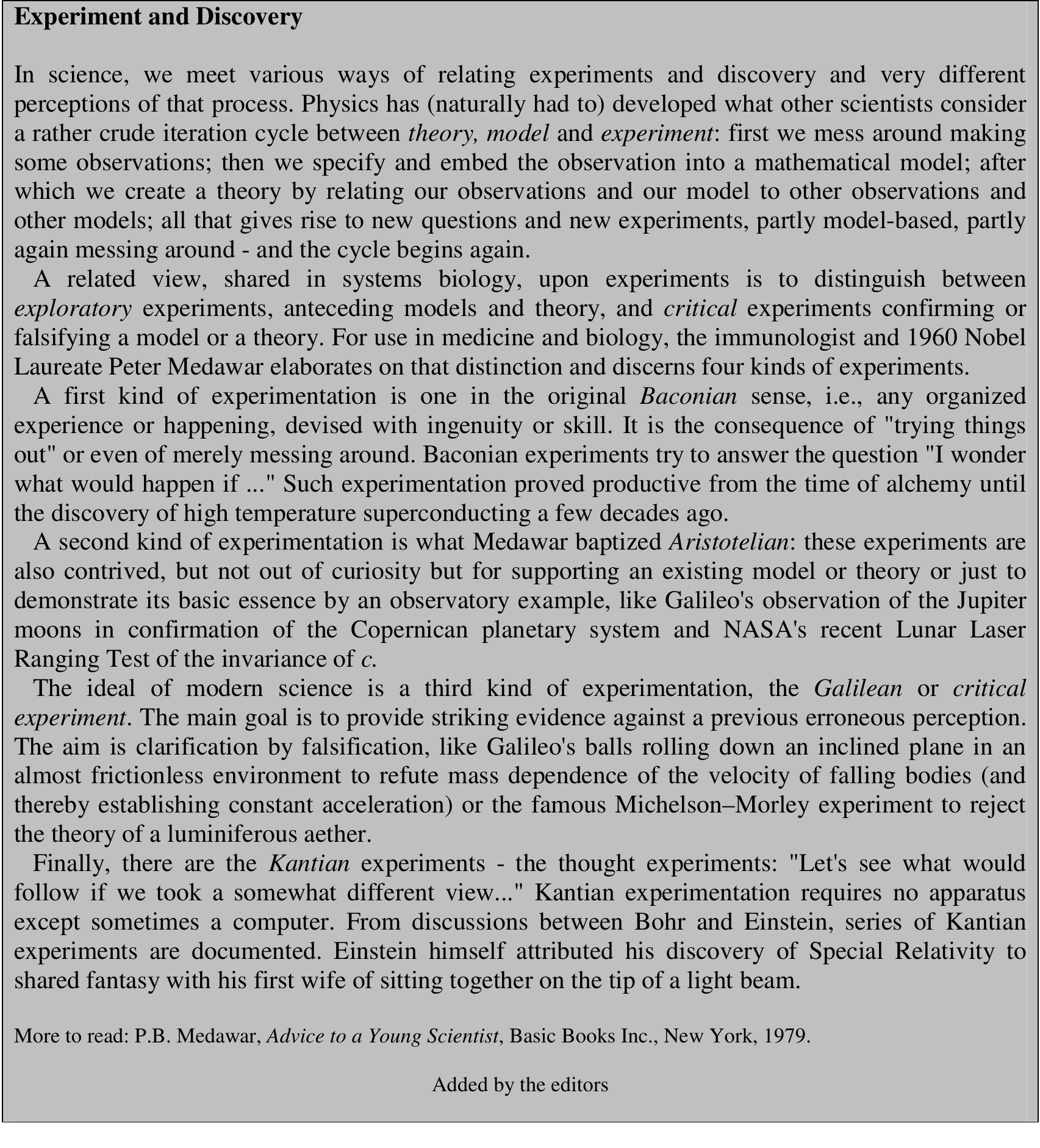}
\end{figure}

\subsection{Suggested Experiments and Measurements}\label{ss:suggested}

At the present stage of our knowledge, the mathematical and
biophysical {\em correctness} of our model does not prove its {\em
relevance} for explaining the phenomena it claims to explain. To
decide whether the here described phenomena and effects are {\em
dispensable} or {\em decisive} for regulated exocytosis, the scales of the
ion oscillations, the electromagnetic fields, the acting forces, the
entering material constants and the characteristic times and lengths
must be determined.

Hence, the framework of our electromagnetic free boundary problem
for the dimple making suggest the following array of experiments and
measurements:\comment{Energy estimates {\`a} la Rosenheck missing}

\begin{enumerate}
\item We shall observe the $\Ca^{2+}$ oscillations prior to the bilayer membrane vesicle fusion also in
pancreatic $\beta$-cells and determine their spatial and temporal character. In particular, the observations
must
\begin{itemize}

\item check the intracellular origin of the oscillations; we may, e.g., deliberately silence
(empty) some types of organelles by adding suitable agents, see, e.g., Fridlyand et al. \cite{FTP};

\item identify the participating organelles ($\Ca^{2+}$ storages); and

\item decide about the orientation (the direction) of oscillations, the frequency, depending on
stimulus, and the distinction between almost simultaneous oscillations
pointing in different directions as their selected sites for the making of the fusion pore.

\end{itemize}

\item We shall modulate the oscillations by submitting the cells to an external field generator with variable
frequencies, to prove the magnetic character of the field wave associated to the oscillations.

\item We shall determine the surface tension in the plasma membrane of living
cells. In particular, we shall measure and/or calculate the bending
rigidity and stretching elasticity under "repacking" of the ball
shaped heads of the inner lipids under area changes. For living
cells, we expect that these magnitudes are substantially larger than
for model membranes, e.g., due to osmotic pressure in living cells,
see Henriksen and Ipsen \cite{HenIps}.

\item We shall measure the cytosol viscosity close to the plasma membrane, i.e.,
update the classic study \cite{BPSV} by Bicknese et al., and locate
actin filaments blocking for unwanted docking of the insulin
granules at the plasma membrane.

\item We shall check Rosenheck's estimate (l.c.) for charged molecules in the plasma membrane.

\item We shall estimate the content of (para-, not necessarily ferro-)magnetizable
$\operatorname{Fe}$ atoms and crystals in the membranes to determine
their magnetic momentum.

\item We shall estimate the distribution of the inhomogeneities of
the magnetic field near the plasma membrane.

\item We need precise electron or atomic force microscope slices of the dimple making and of the degree of its singularity
in $\beta$-cells.

\item We shall measure by patch clamp technique the expected decrease of the capacitive reactance $X_C$
under dimple forming.

\item We shall correlate the $\Ca^{2+}$ oscillations with the fusion events; in particular, we shall
confirm the spatial and temporal coincidence of flickering of exocytosis with break-downs of the field wave.

\item We shall demonstrate the absence or weakness of the $\Ca^{2+}$ oscillations after stimulation
in stressed or tired $\beta$-cells.

\item After obtaining reliable values of all data involved in our
mathematical model, we shall create a computer simulation of the
free boundary problem to get a solution graphically (in form of a
surface). After that, we can compare the dimple images and the
modeled surface.

\end{enumerate}

\begin{acknowledgement}
The first author was partially supported by the Russian Foundation for
Basic Research (grant no. 09-01-00729).
The third author acknowledges the support by the Danish network
{\em Modeling, Estimation and Control of Biotechnological Systems} ({\em MECOBS}).
We four thank the referees for their thoughtful comments to and harsh criticism
of a first draft and F. Pociot and J. St{\o}rling for their corrections and helpful suggestions
which led to many improvements. Referees and colleagues went clearly beyond the call of duty,
and we are indebted to them.

\end{acknowledgement}

\end{document}